\documentstyle[11pt,leqno]{article}
\begin{document}
\newtheorem{theorem}{Theorem}
\newtheorem{lemma}{Lemma}
\newtheorem{corollary}{Corollary}
\newtheorem{proposition}{Proposition}
\newtheorem{definition}{Definition}
\newcommand{\bn}{{\bf n}}
\newcommand{\bT}{{\bf T}}
\newcommand{\noi}{\noindent}
\newcommand{\wtd}{\widetilde}
\newcommand{\wht}{\widehat}
\newcommand{\bm}{{\bf m}}
%\baselineskip = 15pt
%\baseline-stretch = 1.3
\setcounter{page}{0}
\thispagestyle{empty}
\begin{center}
{\bf OPTIMAL STRATEGIES FOR A CLASS OF SEQUENTIAL CONTROL PROBLEMS WITH PRECEDENCE RELATIONS}
\end{center}
\vspace*{0.5cm} \centerline{BY HOCK PENG CHAN$^{*}$, CHENG-DER FUH$^{**}$ AND INCHI
HU$^{***}$} \centerline{\it National University of Singapore, Academia Sinica}
\centerline{\it and Hong Kong University
of Science and Technology} \vspace{1cm} \centerline{\bf Abstract}
\vspace*{0.5cm} 

Consider the following multi-phase project management problem. Each project is 
divided into several phases. All projects enter the next phase at the same point 
chosen by the decision maker based on observations up to that point.  
Within each phase, one can pursue the projects in any order. When pursuing the 
project with one unit of resource, the project state changes according to a 
Markov chain. The probability distribution of the Markov chain is known up to
an unknown parameter. When pursued, the project generates a random reward 
depending on the phase and the state of the project and the unknown parameter.  
The decision maker faces two problems: (a) how to allocate resources to projects 
within each phase, and (b) when to enter the next phase, so that the total 
expected reward is as large as possible. In this paper, we formulate the 
preceding problem as a stochastic scheduling problem and propose asymptotic 
optimal strategies, which minimize the shortfall from perfect information payoff. 
Concrete examples are given to illustrate our method. 

\vspace*{2cm}
\noindent
{\it AMS} 2000 {\it subject classifications.} Primary 62L05;
Secondary 62N99.

\noindent
$^*$ Research supported by grants from the National University of Singapore. \\
$^{**}$ Research partially supported by the National Science Council of ROC.\\
$^{***}$ Research partially supported by Hong Kong Research Grant Council.

\newpage
\def\theequation{1.\arabic{equation}}
\setcounter{equation}{0}
\section{Introduction}
We first formulate the multi-phase project management problem as that of 
optimally scheduling a number of jobs. Suppose that a single machine is 
available to process ${\cal U}$ jobs. Each job belongs to one job group and there 
are $I$ job groups all together. Within each group, the job can be processed in any 
order. However, there exists a predetermined order among job groups. That is, after 
leaving the current job group, there is no return to it in the future.  The state 
of a job under processing evolves as a Markov chain and earns rewards as it is processed, 
not otherwise. The time-varying reward distributions depends on an unknown parameter 
$\theta$. The objective is to minimize the shortfall from perfect information payoff, 
which is the difference between the optimal reward when the parameter is known and 
that when it is unknown. We establish an asymptotic lower bound on this difference 
and construct policies which attain the lower bound. Clearly the preceding stochastic
scheduling problem is the same as the multi-phase project management problem when 
we identify jobs in the same group with projects in the same phase.  

To solve the proposed stochastic scheduling problem, we need to resolve two issues. 
First, our solution must prescribe how to process jobs within the same group. 
Secondly, the solution needs to stipulate the timing of leaving the current job 
group and entering the next one. All existing methods address only one of the two 
issues. As one shall see, to address these two issues simultaneously requires new 
ideas as well as nontrivial combination of existing methods. 

The advantages of efficient strategies constructed in Section 4 is three-fold. 
\begin{itemize}
\item It addresses the two crucial issues described in the previous paragraph 
simultaneously. 
\item It is still optimal, if we consider constant switching cost from one project to another.  
\item When the bad set (see Section 2.4 for definition) is empty 
the strategy is super efficient in the sense of attaining $o(\log N)$ regret 
(see Section 2.2 for definition). 
\end{itemize}

If the parameter $\theta$ were known, the best policy would be to 
process only the job with greatest one-step expected reward. In
ignorance of $\theta$, an optimal policy needs to trade off a reduced reward 
in exchange for information on $\theta$. The key to the optimal trade-off is 
the construction of a strategy that achieves the asymptotic lower bound for 
the shortfall from the complete information payoff, which we shall refer to
as regret hereafter. Although dynamic programming and the Gittins index 
rule (cf. Gittins, 1989) have been developed to solve a general class of 
adaptive control problems, to which the proposed problem belongs, computational 
difficulty makes them less applicable. One reason for adopting the approach 
described here is to obtain an explicit solution which is easy to implement. 

This approach was first introduced by Lai and Robbins (1985) and
generalized by Anantharam, Varaiya and Walrand (1987) and Lai
(1987). When there is only one job group and the rewards from each job
are independent and identically distributed (i.i.d.), the preceding
control problem is the classical multi-armed bandit problem; see
Robbins (1952), Berry and Fristedt (1985) and Gittins (1989). When
there is only one job in each group and rewards are i.i.d., it is  
the irreversible multi-arm bandit problem studied by Hu and Wei (1989),
whereas Hu and Lee (2003) considered the same problem under a Bayesian 
setting. Fuh and Hu (2000) investigated the irreversible multi-armed 
bandit problem with Markovian rewarding. Agrawal, Teneketzis 
and Anantharam (1989a,b) studied controlled i.i.d. processes and Markov chains in  
finite parameter and state spaces. They introduced the concept of bad sets and 
showed that it plays an important role to the solution of the adaptive control
problem. Other related works can be found in Kadane and Simon
(1977), Mandelbaum and Vanderbei (1981), Gittins (1989), Presman and Sonin (1990),
Glazebrook (1991, 1996), Graves and Lai (1997) and references therein.

The rest of the paper is organized as follows. In Section 2, we
describe the components of a statistical model for the proposed problem. 
The asymptotic lower bound for the regret is derived in Section 3. 
In Section 4, we propose a class of strategies making use of an adjusted 
MLE $\widehat \theta_a$. This adjustment is necessary for consistent estimation 
of the bad sets of $\theta$ when the parameter space is continuous. 
The efficiency of our procedure relies on an initial experimentation 
stage based on the adjusted MLE estimate to maximize the information 
content and also on a subsequent testing stage via sequential likelihood 
ratio tests to reject suboptimal jobs or a whole group of jobs. 
Unequal allocation of processing time on jobs may occur in the testing stage 
so that there is more frequent processing of superior jobs. In Section 5, 
we discuss how our method can be applied to multi-phase project management
examples. Most of the technical proofs are deferred to the Appendix.

\def\theequation{2.\arabic{equation}}
\setcounter{equation}{0}
\section{Preliminaries}
\subsection{The scheduling problem}

Let ${\cal U} = J_1 + \cdots +J_I$ indicate that there are $I$ groups and 
$J_i$ jobs in the $i$th group for $i=1,\ldots,I$. One is free to process any job 
within the same group, while jobs must be processed following the order of
$1,\ldots,I$ between groups. As processing a job a unit time is equivalent to
taking an observation from a statistical population, we have
${\cal U}$ statistical populations $\Pi_{11},\ldots,\Pi_{IJ_I}$. For each
$ij$, the observations from $\Pi_{ij}$ follow a Markov chain on a state 
space $D$ with $\sigma$-algebra ${\cal D}$. It is assumed that the transition
probability $P_{ij}^{\theta}$ for the Markov chain has a probability
density function $p_{ij}(x,y;\theta)$ with respect to some
nondegenerate measure $Q$, where $p_{ij}(x,y;\cdot)$ is known and
$\theta$ is an unknown parameter belonging to a parameter space
$\Theta$. We assume that the stationary probability
distribution for the Markov chain exists and has probability density
function $\pi_{ij}(\cdot;\theta)$ with respect to $Q$. 
%Let $g_{ij}(x)$ be the reward when population $\Pi_{ij}$ is sampled and $x$ is observed. 
At each step, we are required to process one job respecting the partial 
order $ij\preceq i'j'\Leftrightarrow i\leq i'$. 

An adaptive policy is a rule that dictates, at each step, which job should be
processed based on information from previous observations. We can represent
a policy as a sequence of random variables $\phi = \{ \phi_t \}$ taking values in 
$\{ij:i=1,\cdots,I;~ j=1,\cdots, J_i\}$, such that the event $\{\phi_t = ij \}$
(process job $ij$ at step $t$) belongs to the $\sigma$-field ${\cal F}_{t-1}$ generated by 
$\phi_1,X_1,\ldots,\phi_{t-1}, X_{t-1}$, where $X_n$ denotes the
state of the job being processed at the $n$th step. The constraint
\begin{eqnarray}\label{poc}
   \phi_t \preceq \phi_{t+1}~~~~~~\mbox{for}~~1 \leq t \leq N-1,
\end{eqnarray}
indicates that once a sample has been taken from $\Pi_{ij}$, one can
switch to other jobs within group $i$ or to the jobs in groups
$i+1$ to $I$, but no further sampling is allowed from
$\Pi_{11},\ldots,\Pi_{(i-1)J_{i-1}}$. 

%Note that the jobs considered here are ``{\it correlated}'', that is, 
%when a job is processed, in addition to providing information on the reward
%distribution of the processed job, it also gives distributional
%information about other jobs as well. The difference between the
%``{\it correlated}'' and the ``{\it independent}'' jobs lies in
%the structure of the parameter space $\Theta$. In the correlated-
%job problem, the parameter $\theta \in \Theta$
%parameterizes $all$ jobs $\Pi_{ij}$, whereas in the independent
%job problem, the parameter $(\theta_{11},\cdots,\theta_{IJ_I}) \in \Theta$ 
%is such that each $\theta_{ij}$ parameterizes the $individual$ $\Pi_{ij}$.

\subsection{The objective function}
Let the initial state of the job $ij$ under processing be distributed
according to $\nu_{ij}(\cdot;\theta)$. Throughout this paper, we shall use
the notation $E_\theta$ ($P_\theta$) to denote expectation (probability)
with respect to the initial distribution $\nu_{ij}(\cdot;\theta)$; 
similarly, $E_{\pi(\theta)}$ to denote expectation with respect to $P_\theta$
and the stationary
distribution $\pi_{ij}(\cdot;\theta)$. We shall assume that ${\cal V}_{ij} = \{
x \in D: \nu_{ij}(x;\theta) > 0 \}$ does not depend on $\theta$ and
\begin{equation} \label{2.2a}
v_{ij}:= \inf_{x \in {\cal V}_{ij}} \inf_{\theta,\theta' \in \Theta} [\nu_{ij}
(x;\theta)/\nu_{ij}(x;\theta')] > 0 \mbox{ for all } i,j. 
\end{equation}
Suppose that
$\int_{x \in D} |g(x)| \pi_{ij}(x;\theta) Q(dx) < \infty$ for some real-valued
function (reward) $g$. Let
\[
\mu_{ij} (\theta) = \int_{x \in D} g(x) \pi_{ij}(x;\theta) Q(dx)
\]
be the mean reward under stationary distribution $\pi_{ij}$ if job $ij$ 
is processed once. Let $N$ be the total processing time for all jobs, and
\begin{eqnarray}\label{TNij}
T_N(ij) = \sum_{t=1}^N {\bf 1}_{\{\phi_t = ij \}}
\end{eqnarray}
be the amount of time that job $ij$ is processed and ${\bf 1}$ denotes
the indicator function. An optimal strategy would be one which maximizes
\begin{eqnarray}\label{totalreward} \qquad 
 W_N (\theta) :=
 \sum_{t=1}^N \sum_{i=1}^I \sum_{j=1}^{J_i} E_{\theta} \{ E_{\theta}
 [g(X_t) {\bf 1}_{\{\phi_t = ij\}} | {\cal F}_{t-1} ] \}.
\end{eqnarray}
In the case of independent rewards, that is, when $p_{ij}(x,y,;\theta)
= p_{ij}(y;\theta)$ for all $i,j,x,y$ and $\theta$, $W_N(\theta) 
= \sum_{i=1}^I \sum_{j=1}^{J_i} \mu_{ij}(\theta) E_{\theta} T_N(ij)$. We shall
show in the Appendix that for Markovian rewards, under regularity conditions
A3-A4 (see Section 2.3), there exists a constant $C_0 < \infty$ independent of $\theta \in
\Theta$, $N > 0$ and the strategy $\phi$ such that
\begin{equation} \label{2.3a}
\Big| W_N(\theta) - \sum_{i=1}^I \sum_{j=1}^{J_i} \mu_{ij}(\theta) E_{\theta} T_N(ij)
\Big| \leq C_0.
\end{equation}
When the parameter space $\Theta$ and state space $D$ are both finite,
(\ref{2.3a}) also follows from Anantharam, Varaiya and Walrand (1987, Lemma 2.1). 
In light of (\ref{2.3a}), maximizing $W_N(\theta)$ is asymptotically
equivalent [up to a $O(1)$ term] to minimizing the regret
\begin{eqnarray}\label{regret} 
 R_N (\theta) & := & N \mu^*(\theta) -  
\sum_{i=1}^I \sum_{j=1}^{J_i} \mu_{ij}(\theta) E_{\theta} T_N(ij) \cr
& = & \sum_{ij:\mu_{ij}(\theta) < \mu^*(\theta)}
   [\mu^*(\theta) - \mu_{ij}(\theta)] E_{\theta} T_N(ij),
\end{eqnarray}
where $\mu^*(\theta) := \max_{1 \leq i \leq I} \max_{1 \leq j \leq J_i}
\mu_{ij}(\theta)$.

Because adaptive strategies $\phi$ which are optimal for all
$\theta \in \Theta$ and large $N$ in general do not exist, we
consider the class of all (asymptotically) {\it uniformly good}
adaptive strategies under the partial order constraint $\preceq$,
with regret satisfying
\begin{equation}\label{uniformlygood}
R_N (\theta) = o(N^{\alpha}), \quad \mbox{ for all } \alpha > 0 \mbox{ and } 
\theta \in \Theta.
\end{equation}
Such strategies have regret that does not increase too rapidly for
any $\theta \in \Theta$. We would like to find a strategy
that minimizes the increasing rate of the regret within the
class of uniformly good adaptive strategies under the partial order
constraint $\preceq$. 

Due to the irreversibility constraint (\ref{poc}),
a strategy satisfying (\ref{uniformlygood}) would in general be dependent on $N$
when there are more than one group of arms. Consider for example the case in which
the optimal arm is unique and lies in the first group. Let $p > 0$ be
the probability that a strategy $\phi$ bypasses the first group of arms
before a fixed time $N_0$. If the strategy $\phi$ is independent of $N$, then
$$R_N(\theta) \geq p(N-N_0)[\mu^*(\theta) - \max_{2 \leq i \leq I}
\max_{1 \leq j \leq J_i} \mu_{ij}(\theta)]
$$
and (\ref{uniformlygood}) does not hold. This is unlike the case of  
one group multi-armed bandit considered by Lai and Robbins (1985), Anantharam et
al. (1987) and Agrawal et al. (1989a,b) whereby optimal strategies $\phi$ satisfying 
(\ref{uniformlygood}) and not dependent on $N$ have been constructed.

%Let $a(\theta) > 0$ be the switching cost between two 
%arms which are not both optimal when the underlying parameter is $\theta$. It
%is assumed here that there is no switching cost when both arms are optimal.
%Then
%$$H_N(\theta) := a(\theta) E_\theta \Big( \sum_{t=1}^{N-1} {\bf 1}_{\{
%\phi_t \neq \phi_{t+1}, \min[ \mu_{\phi_t}(\theta), \mu_{\phi_{t+1}}(\theta) ]
%< \mu^*(\theta) \}} \Big)$$
%is the average switching cost of a procedure. It is also
%desirable that this cost is asymptotically negligible compared to the regret
%as $N \rightarrow \infty$.
  
\subsection{The assumptions}
Denote the Kullback-Leibler information number by
\begin{eqnarray}\label{definekl}
 I_{ij}(\theta,\theta') =
   \int_{x \in D} \int_{y \in D}  \log \Big[ \frac{p_{ij}(x,y;\theta)}
{p_{ij}(x,y;\theta')} \Big] p_{ij}(x,y;\theta)\pi_{ij}(x;\theta) Q(dy) Q(dx).
\end{eqnarray}
%where $p_{ij}(x,y;\theta)$ is the transition density from state $x$ to state $y$ 
%of the $ij$th Markov chain. 
Then, $0 \leq I_{ij}(\theta,\theta') \leq \infty$. We shall
assume  that $I_{ij}(\theta,\theta') < \infty$ for all $i,j$ and
$\theta,\theta' \in \Theta.$ Let $\mu_i(\theta) = \max_{1 \leq j \leq J_i}
\mu_{ij}(\theta)$ be the largest reward in the $i$th group of jobs, and 
\begin{equation}\label{Thetai}
 \Theta_i = \{ \theta \in \Theta: \mu_i (\theta)> \mu_{i'}
 (\theta)~{\rm for \ all}~i'<i~{\rm and} \ \mu_i(\theta) \geq \mu_{i'}
(\theta)~{\rm for \ all}~i'\geq i \}
\end{equation}
be the set of parameter values such that the first optimal job is in group $i$. Let 

\begin{equation}
\Theta_{ij} = \{ \theta \in \Theta_i: \mu_{ij}(\theta) = \mu_i(\theta)
\}
\end{equation}
be the parameter set such that job $ij$ is one of the first optimal jobs.  
Each $\theta \in \Theta$ belongs to exactly one $\Theta_i$ but may belong
to more than one $\Theta_{ij}$. Let
\begin{equation}
\Theta_i^* = \{ \theta \in \Theta: \mu_i(\theta) > \mu_{i'}(\theta) \ {\rm for \ all
} \ i' \neq i \}
\end{equation}
be the parameter set in which all the optimal arms lie in group $i$.
Clearly, $\Theta_i^* \subset \Theta_i$ but the reverse relation is not
necessarily true.

We now state a set of assumptions that will be used to prove the optimality 
results in Sections 3 and 4. Let $\Theta$ be a compact subset of ${\bf R}^d$ for some $d \geq 1$
and let $X_{ijt}$ denotes the $t$th observation taken from arm $ij$. 

\begin{itemize}
\item [A1.] $\mu_{ij}(\cdot)$ are finite and continuous on $\Theta$ for all $i,j$. 
Moreover, no job group is redundant in the sense that $\Theta^*_i\neq\emptyset$ for all $i=1,\cdots,I$. 
%$\Theta_i^*$ is a dense subset of $\Theta_i$ for all $i$. In other words,
%for every $\theta_i \in \Theta_i$, there exists $\theta_i^{(k)} \in
%\Theta_i^*$ for $k=1,2,\ldots$ such that $\lim_{k \rightarrow \infty}
%\theta_i^{(k)} = \theta_i$.

\item [A2.] $\sum_{j=1}^{J_1}
I_{1j}(\theta,\theta') > 0$ for all $\theta' \neq \theta$ and
$\inf_{\theta' \in \Theta_{ij}} I_{ij}(\theta,\theta') > 0$ for all  
$1\leq i<I,1\leq j\leq J_i$ and $\theta \in \cup_{\ell>i} \Theta_\ell$.

\item [A3.] For each $j=1,\ldots,J_i,i=1,\ldots,I$ and $\theta \in
\Theta$, $\{X_{ijt},t \geq 0\}$ is a Markov chain on a state
space $D$ with $\sigma$-algebra ${\cal D}$, irreducible
with respect to a maximal irreducible measure on $(D, {\cal D})$ and
aperiodic. Furthermore, $X_{ijt}$ is {\it Harris
recurrent} in the sense that there exists a set $G_{ij}\in {\cal
D}$, $\alpha_{ij} > 0$ and probability
measure $\varphi_{ij}$ such that $P_{ij}^{\theta}\{
X_{ijt} \in G_{ij}~ {\rm i.o.}|X_{ij0}=x\}=1$ for all $x \in D$ and
\begin{eqnarray}\label{hrecurrent} \qquad
P^{\theta}_{ij} \{ X_{ij1} \in A|X_{ij0}=x \} \geq \alpha_{ij}
\varphi_{ij}(A) \quad {\rm for \ all} \ x \in G_{ij} \ {\rm and} \ A \in {\cal D}.
\end{eqnarray}

\item [A4.] There exist constants $0 < \bar{b} < 1$, $b > 0$ and drift functions $V_{ij}: D \rightarrow
[1,\infty)$ such that for all $i=1,\ldots,I$ and $j=1,\ldots,J_i$, 
\begin{equation} \label{2.10a}
\sup_{x \in D} |g(x)|/V_{ij}(x) < \infty,
\end{equation}
and for all $x \in D$ and $\theta \in \Theta$,
\begin{equation}\label{drift}
P_{ij}^{\theta}V_{ij}(x) \leq (1-\bar{b}) V_{ij}(x) +b {\bf 1}_{\{x \in G_{ij} \}},
\end{equation}
where $G_{ij}$ satisfies
(\ref{hrecurrent}) and $P_{ij}^{\theta}V_{ij}(x)=\int_D V_{ij}(y) P_{ij}^{\theta}(x,dy)$.
Moreover, we require that 
\begin{equation}\label{vij}
\int_D V_{ij}(x) \nu_{ij}(x;\theta) Q(dx) <\infty \quad {\rm and} \ V^*_{ij} :=
\sup_{x \in G_{ij}} V_{ij}(x) < \infty.
\end{equation}

Let $\ell_{ij}(x,y;\theta,\theta') = \log [p_{ij}(x,y;\theta)/p_{ij}(x,y;\theta')]$
be the log likelihood ratio between $P^\theta_{ij}$ and $P^{\theta'}_{ij}$ 
and $N_\delta(\theta) = \{ \theta': \| \theta - \theta' \| < \delta \}$ 
a ball of radius $\delta$
around $\theta$, where $\| \cdot \|$ denotes Euclidean norm.

\item [A5.] There exists $\delta > 0$ such that for all $\theta, \theta' 
\in \Theta$, 
\begin{equation}\label{a51}
K_{\theta,\theta'} := 
\sup_{x \in D} {E_{\theta} [ \sup_{\tilde \theta \in N_{\delta}(\theta')}
\ell_{ij}^2(X_{ij0},X_{ij1};\theta,\wtd \theta)|X_{ij0}=x] \over V_{ij}(x)} < \infty
\end{equation}
for all $j=1,\ldots,J_i$, $i=1,\ldots,I$. Moreover, 
\begin{equation}\label{a52}
\sup_{\tilde \theta \in N_{\delta'}(\theta')} | \ell_{ij}(x,y;\theta',\wtd \theta) |
\rightarrow 0 \ {\rm as} \  \delta' \rightarrow 0
\end{equation}
for all $x,y \in D$ and $\theta' \in \Theta$. 

\end{itemize}

Assumption A1 is for excluding some unrealistic models in which efficient but
impractical strategies may exist. A2 is a positive information criterion: the first inequality 
makes sure that information is available in the first job group to estimate $\theta$; 
while the second inequality allows us to gather information in the $i$th job group 
for moving to the next group when $\theta \in \Theta_\ell$ for some $\ell > i$. 
Assumption A3 is a recurrence condition and A4 is a drift condition. These two conditions 
are used to guarantee the stability of the Markov chain so that the strong law 
of large numbers and Wald's equation hold. A5 is a finite second moment condition that
allows us to bound the probability that the MLE of $\theta$ lies outside a 
small neighborhood of $\theta$. This bound is important for us to determine the level 
of unequal allocation of observations that can be permitted in the testing stage of 
our procedure. The proof of the asymptotic lower bound in Theorem 1 requires only A1-A3; 
while additional A4 and A5 are required for the construction of efficient strategies 
attaining the lower bound.

We now demonstrate an immediate consequence of A3-A5 that for any $\theta \in \Theta$ and
$\varepsilon > 0$, there exists $0 < \delta' < \delta$ such that
\begin{equation}\label{localsup}
E_{\pi_{ij}(\theta)} \Big[ \sup_{\tilde \theta \in N_{\delta'}(\theta')}
| \ell_{ij}(X_{ij0},X_{ij1};\theta',\wtd \theta) | \Big] < \varepsilon
\end{equation}
for all $ij$ and $\theta' \in \Theta$. Note that the continuity of $I_{ij}(\theta,\cdot)$
follows from (\ref{localsup}). 

Since $\pi_{ij} = C' \sum_{k=0}^\infty (P_{ij}-\alpha_{ij}\varphi_{ij} {\bf 1}_{G_{ij}})^k 
\varphi_{ij}$, where $C'$ is a normalizing constant, it follows from (\ref{drift})-(\ref{vij}) 
that $\int_D V_{ij}(x) \pi_{ij}(x;\theta) Q(dx)
< \infty$. Hence by (\ref{a51}) and the relation $\ell_{ij}(X_{ij0},X_{ij1};\theta',\wtd \theta)
= \ell_{ij}(X_{ij0},X_{ij1};\theta,\wtd \theta)-\ell_{ij}(X_{ij0},X_{ij1}
;\theta,\theta')$, we have
\begin{eqnarray*}
& & E_{\pi_{ij}(\theta)} \Big[ \sup_{\tilde \theta \in N_{\delta}(\theta')}
| \ell_{ij}(X_{ij0},X_{ij1};\theta',\wtd \theta) | \Big] \cr
& & \quad \leq 
E_{\pi_{ij}(\theta)} |\ell_{ij}(X_{ij0},X_{ij1};\theta,\theta')| +  
E_{\pi_{ij}(\theta)} \Big[ \sup_{\tilde \theta \in N_\delta(\theta')} |\ell_{ij}(X_{ij0},
X_{ij1};
\theta,\wtd \theta)| \Big] < \infty.
\end{eqnarray*}
As the convergence in (\ref{a52}) is monotone decreasing, 
it follows from the dominated convergence theorem that (\ref{localsup}) holds. 

\subsection{Bad sets}
Bad set is a useful concept for understanding the learning required within the group 
containing optimal jobs. It is associated with the asymptotic lower bound described in Section 3 
and is used explicitly in Section 4 to construct the asymptotically efficient strategy.
For $\theta \in \Theta_\ell$, define 
$J(\theta) = \{ j: \mu^*(\theta) = \mu_{\ell j}(\theta) \}$ as the set of
optimal jobs in group $\ell$. Hence $\theta \in \Theta_{\ell j}$ if and only if
$j \in J(\theta)$. We also define the bad set, the set of `bad'
parameter values associated with $\theta$, as all $\theta' \in \Theta_\ell$ which
cannot be distinguished from $\theta$ by processing any of the optimal jobs $\ell j$.
More specifically, the bad set
\begin{equation}\label{bs}
 B_\ell(\theta) = \Big\{ \theta' \in \Theta_\ell \setminus
 \big( \bigcup_{j \in J(\theta)} \Theta_{\ell j} \big) :~I_{\ell 
j}(\theta,\theta') = 0 \ {\rm for \
all} \ j \in J(\theta) \Big\}. 
\end{equation}
We note that if $I_{\ell j}(\theta,\theta')=0$, then the transition probabilities 
of $X_{\ell jt}$
are identical under both $\theta$ and $\theta'$. If $\theta' \in
B_\ell(\theta)$, then by definition, $\theta' \not\in \cup_{j \in J(\theta)}
\Theta_{\ell j}$ and hence $J(\theta') \cap J(\theta) = \emptyset$.
Let $j \in J(\theta)$ and $j' \in J(\theta')$. Then
$\mu_{\ell j'}(\theta') > \mu_{\ell j}(\theta') = \mu_{\ell j}(\theta) >
\mu_{\ell j'}(\theta)$. Thus
\begin{equation}\label{bsdual}
I_{\ell j'}(\theta, \theta') > 0 \ {\rm for \ all} \ \theta' \in B_\ell(\theta) \ 
{\rm and} \ j' \in J(\theta').
\end{equation}
The interpretation of (\ref{bsdual}) is as follows. Although we cannot distinguish $\theta$ from 
$\theta'\in B_{\ell}(\theta)$ when processing the optimal job for $\theta$, we can distinguish them 
by processing the optimal job for $\theta'$. 
This fact explains the necessity of 
processing non-optimal jobs to collect information. 

Assumption A2 says when sampling from the optimal arm one {\it can} distinguish
any $\theta$ value whose optimal arm is in a {\it future} group. But having
a non-empty bad set says that when sampling from the optimal arm 
{\it cannot} distinguish some $\theta$ value whose optimal
arm is in the {\it current} group. These two statements are compatible.
We now provide two examples from 
the celebrated multi-armed bandit problem to illustrate the idea of bad sets. \\

\noindent {\bf Example 1: Independent armed-bandit problem.} 
Let $\Pi_{11},\ldots,\Pi_{1J}$ denote $J$ statistical
populations specified, respectively, by density functions
$p(x;\theta_j)$ with respect to some measure $Q$. For simplicity,
assume that $x=0,1$ and  $p(0;\theta_j) = 1-\theta_j$, $p(1;\theta_j) =
\theta_j$, where $\theta_j$ are unknown parameters taking values
in $[0,1]$. A multi-armed bandit problem searches for strategies to
sample $X_1,X_2,\ldots,$ sequentially from these $J$ populations in order
to maximize the expected value of the sum $S_N = \sum_{t=1}^N X_t$
as $N \rightarrow \infty$.

Let $\theta = (\theta_1,\ldots, \theta_J)$.
If $\theta=(0.2,0.1)$, then the set of optimal arms 
$J(\theta)= \{ 1 \}$ and the bad set $B_1(\theta)=\{ (0.2,\theta_2'): 0.2 < 
\theta_2' \leq 1 \}$. Even though arm 1 is optimal, experimentation from
arm 2 is required to make sure that the true parameter value does not lie in
$B_1(\theta)$.

The two-armed bandit problem studied by Feldman (1962) has $\Theta=\{(\theta_1, 
\theta_2)$, $(\theta_2,\theta_1)\}$ with $\theta_1\neq \theta_2$. It follows that  
$B_1(\theta) = \emptyset$ for all 
$\theta \in \Theta$. This leads to remarkably low regret, $R_N(\theta) =O(1)$.
\\

\noindent
{\bf Example 2: Correlated armed-bandit problem.}
Consider bivariate normal populations
$\Pi_{11},\Pi_{12},\Pi_{13}$ with respective mean vectors
$(\mu_1,\lambda),$ $(\mu_2,\mu_3)$ and $(\mu_3, \mu_2+\lambda)$,
where $\mu_1,\mu_2,\mu_3,\lambda$ are unknown parameters. The problem is
to sample the random vectors sequentially to maximize the
expected value of the first component of the observed sum,
$\sum_{t=1}^N X_t$, as $N \rightarrow \infty$. Let
$\theta = (\mu_1,\mu_2,\mu_3,\lambda)$. If
$J(\theta) = \{ 1 \}$, then 
$$B_1(\theta) = \{ \theta' \in \Theta: \mu_1 = \mu_1',
\ \lambda = \lambda', \ \max(\mu_2',\mu_3') > \mu_1' \}.$$

\section{A lower bound for the regret}
\def\theequation{3.\arabic{equation}}
\setcounter{equation}{0}
The following theorem gives an asymptotic lower bound for the regret (\ref{regret})
of uniformly good adaptive strategies under the partial order constraint
$\preceq$.

\begin{theorem}
Assume {\rm A1}-{\rm A3} and let $\theta \in \Theta_\ell$.
For any uniformly good adaptive strategy
$\phi$ under the  partial order constraint $\preceq$, 
\begin{eqnarray}\label{lb}
\liminf_{N \rightarrow \infty} {R_N(\theta)}/{\log N} \geq  z(\theta,\ell),
\end{eqnarray}
where $z(\theta,\ell)$ is a solution of the following minimization problem. 
\begin{equation}\label{mp}
\mbox{Minimize } \sum_{i<\ell} \sum_{j=1}^{J_i}[\mu^*(\theta)-
    \mu_{ij} (\theta)] z_{ij}(\theta) + \sum_{j\notin J(\theta)}
[\mu^{*}(\theta)-\mu_{\ell j}(\theta)] z_{\ell j}(\theta),
\end{equation}
subject to $z_{ij}(\theta)\geq 0$, 
$j = 1,\ldots,J_{i}$ if $i < \ell$; $j \notin J(\theta)$ if $i=\ell$; and
 \begin{equation} \label{mpconstraint} \qquad
 \left\{ \begin{array}{lll}
      \inf_{\theta' \in \Theta_1}
       \{ \sum_{j=1}^{J_1} I_{1j}(\theta,\theta')z_{1j}(\theta)\} \geq 1, \\
      \inf_{\theta' \in \Theta_2}
       \{ \sum_{j=1}^{J_1} I_{1j}(\theta,\theta')z_{1j}(\theta)+
       \sum_{j=1}^{J_2} I_{2j}(\theta,\theta')z_{2j}(\theta) \} \geq 1, \\
     \vdots  \\
      \inf_{\theta' \in \Theta_{\ell-1}}
       \{ \sum_{j=1}^{J_1} I_{1j}(\theta,\theta')z_{1j}(\theta)+ \cdots
       + \sum_{j=1}^{J_{\ell-1}} I_{(\ell-1)j}(\theta,\theta')
z_{(\ell-1)j}(\theta) \} \geq 1,\\
   \inf_{\theta'\in B_\ell(\theta)}
       \{\sum_{i < \ell} \sum_{j=1}^{J_i} I_{ij}(\theta,\theta')z_{ij}(\theta)
+\sum_{j \notin J(\theta)}I_{\ell j}(\theta,\theta')z_{\ell j}
       (\theta) \} \geq 1. \end{array}    \right.   \end{equation}

\end{theorem}

The first $(\ell-1)$ inequalities in (\ref{mpconstraint}) are due to the partial order 
constraints. When there is no partial order constraint and the jobs are 
independent, the solution of Problem A reduces to the lower bound given 
in Theorem 1 of Lai and Robbins (1985).  
%Agrawal, Teneketzis and Anantharam (1989) studied a
%controlled independent process with finite parameter space, and
%introduced a finite set $B(\theta)$ of bad parameter values associated
%with $\theta$. In (2.8), for a general parameter space $\Theta$,
%we define the bad set $B(\theta)$ via the decomposition
%$\Theta=\cup_{i=1}^{I} \cup_{j=1}^{J_i} \Theta_{ij}$ and apply it to
%Theorem 1 for the regret lower bound of bandits under partial order
%constraint.

%For a uniformly good adaptive strategy and $\theta \in \Theta$, we have
%\begin{eqnarray}
%\lim_{N \rightarrow \infty} N^{-1} W_N (\theta) = \mu^*(\theta).
%\end{eqnarray}
%The strategies that satisfy (3.3) are said to be {\em consistent}.
Under the assumptions of Theorem 1, the strategies that satisfy, 
for $\theta \in \Theta_\ell,$
\begin{eqnarray}\label{aef}
% \sum_{i>\ell} \sum_{j=1}^{J_i} E_{\theta} T_N(ij) &=& O(\log N),  \\
\lim_{N \rightarrow \infty} R_N(\theta)/\log N &=& z(\theta,\ell), 
\end{eqnarray}
are said to be {\em asymptotically efficient}. If $B_{\ell}(\theta)=\emptyset$, then
the last inequality of (\ref{mpconstraint}) is removed. In particular, when 
$\theta \in \Theta_1,$ (\ref{aef}) implies that
\begin{eqnarray}\label{bsregret}
 R_N(\theta) = \cases{ O(\log N) & {\rm if} \ $B_1(\theta) \neq \emptyset$, \cr
o(\log N) & {\rm if} \ $B_1(\theta) = \emptyset$.}
\end{eqnarray}
We shall assume that $B_{\ell}(\theta)$ is non-empty for the underlying $\theta\in \Theta_\ell$, 
which is true for most applications. The case of $B_\ell(\theta)=\emptyset$ will be
treated elsewhere. 
%the case of correlated arms considered in Agrawal, TeneketzisAnantharam (1989). 
%It also follows from (\ref{aef}) that
%\begin{eqnarray}
% E_{\theta} \Big[ N - \sum_{j \in J(\theta)} T_N(\ell j) \Big] = O(\log N).
%\end{eqnarray}

\medskip
The following lemma will be used to prove Theorem 1. The proofs of
both Lemma 1 and Theorem 1 will be given in the Appendix. 

\begin{lemma}~Assume {\rm A2}-{\rm A3}.
Let $\phi$ be a uniformly good adaptive strategy under the partial order
constraint $\preceq$. If $\theta \in \Theta_{\ell}$,
then for every $\theta' \in \Theta_k^*$, $k < \ell$,
\begin{eqnarray}\label{ml}
\liminf_{N \rightarrow \infty} \Big\{ \sum_{i=1}^k \sum_{j=1}^{J_i}
I_{ij} (\theta,\theta') E_{\theta} T_N(ij)  \Big\} / \log N \geq 1,
\end{eqnarray}
and for every $\theta' \in B_{\ell}(\theta)$,
\begin{eqnarray}\label{mlbs} \qquad
\liminf_{N \rightarrow \infty} \Big\{\sum_{i< \ell} \sum_{j=1}^{J_i}
I_{ij}(\theta,\theta')E_{\theta} T_N(ij)+\sum_{j \notin J(\theta)}
I_{\ell j} (\theta,\theta') E_{\theta} T_N(\ell j) \Big\}/\log N
\geq 1.
\end{eqnarray}
\end{lemma}

%Applying Lemma 1 for $\ell=1,\ldots,I$, we obtain
%\begin{corollary}
%Assume {\rm A2}-{\rm A3}. Let $\phi$ be an uniformly good
%strategy  under the partial order constraint $\preceq$. Then, for
%every $\theta \in \Theta_\ell$, $\theta_k \in \Theta_k^*$ for $k < \ell$
%and $\theta_\ell \in B_\ell(\theta)$, 

%\[  \left\{ \begin{array}{ll}
%\liminf_{N \rightarrow \infty} \sum_{j=1}^{J_1} I_{1j}(\theta,\theta_1)
%   E_{\theta}T_{N}(1j)  /\log N \ge 1,    \\
%  \vdots   \\
%  \liminf_{N\rightarrow\infty} \sum_{i=1}^{\ell-1} \sum_{j=1}^{J_i}
%  I_{ij}(\theta,\theta_{\ell-1}) E_{\theta}T_{N}(ij)/\log N \ge 1, \\
%\liminf_{N \rightarrow \infty} \Big\{\sum_{i< \ell} \sum_{j=1}^{J_i}
%I_{ij}(\theta,\theta_\ell)E_{\theta} T_N(ij) \\
%~~~~~~~~~~~~~~~~+ \sum_{j \notin J(\theta)}
%I_{\ell j} (\theta,\theta_\ell)E_{\theta}T_N(\ell j) \Big\}/\log N \ge 1.
%\end{array}\right.
%\]
%\end{corollary}

\section{Construction of asymptotically efficient strategies}
\def\theequation{4.\arabic{equation}}
\setcounter{equation}{0}

\subsection{Outline of the construction}
The goal of any reasonable strategy is to determine whether the job
currently under processing is optimal or not based on sequential
observations. The job under processing, say job $ij$, is optimal
if $\theta\in\Theta_{ij}$. Thus, the problem of constructing an
efficient adaptive strategy reduces to that of finding a procedure  
to determine whether $\theta \in \Theta_{ij}$ is true or not based on a
sequential sample. The asymptotic lower bound discussed in Section 3
gives us valuable information about the size of the sequential
sample. In particular, it suggests that for $\theta \in
\Theta_\ell$, the amount of processing time for job $ij$,
$j=1,\ldots,J_i$, $i<\ell$, and $j \notin J(\theta)$ if $i=\ell$ should be
$[z_{ij}(\theta)+ o(1)] \log N$, where $z_{ij}(\theta)$ solves the 
minimization problem (\ref{mp}). 

In view of Theorem 1, the sample size $[z_{ij}(\theta)+ o(1)] \log N$ 
represents the minimum amount of learning about job $ij$ in order for 
the strategy to be uniformly good. Because of the partial order constraint
$\preceq$, we also need a sequential test to ensure that the optimal job 
is passed over with probability not exceeding $N^{-1}$. These two facts 
are important guidelines for the construction of asymptotically efficient 
strategies so that the two crucial issues mentioned in the abstract and 
Section 1 can be addressed. 
 
%Recall that $T_N(ij)$ defined in (2.2) is the amount of processing time
%for job $ij$. Let $T_N(i)=T_N(i1)+\cdots + T_N(iJ_i)$ be the amount
%of processing time for group $i$. Note that $\phi_{t}$ in group $\ell$ 
%implies that
%\begin{equation}
%\sum^{\ell-1}_{i=1} \sum_{j=1}^{J_i}T_N(ij)< t \le \sum^{\ell}_{i=1} \sum_{j=1}^{J_i}
%T_N(ij).
%\end{equation}
%Our goal is to construct stopping-times
%$T_N(ij)$ so that the regret 
%$R_N(\theta) \sim z(\theta,\ell) \log N$ for
%$\theta \in \Theta_\ell$, 
%where $z(\theta,\ell)$ is the solution to (\ref{mp}).

 Let $n_0, n_1$ be positive integers that increase to infinity with respect
to $N$ such that $n_0 = o(\log N)$ and $n_1=o(n_0)$. We shall now describe 
the asymptotically efficient strategy $\phi^*$ by dividing it into 
three distinct stages; estimation, experimentation and testing. 

In the estimation stage, $n_0 =o(\log N)$ observations are taken 
from each job in group 1 for estimating the parameter $\theta \in \Theta_\ell$.
If $\ell > 1$ or $\ell=1$ and $B_1(\theta) \neq \emptyset$, 
then an order of $\log N$ observations are taken in the experimental 
stage which contribute $[z(\theta,\ell)+o(1)] \log N$ to the regret; see (\ref{lb}). 
Finally, in the testing phase, $o(\log N)$ observations are taken from each of 
the suboptimal jobs. We first consider the optimal strategy for the case of finite $\Theta$, 
which captures the essential ingredients without too much technical details.
We then extend the strategy to infinite $\Theta$ followed by a formal statement of 
optimality in Theorem 2.  

\subsection{\bf Optimal strategy for finite $\Theta$}

1. {\it Estimation}. Take an initial sample of $n_0$ observations from each job in group 1.  
Let $\widehat{\theta}$ be the maximum likelihood estimate (MLE) of 
$\theta$ defined by 
\begin{equation}\label{mle}
L(\theta)  = \sum_{j=1}^{J_1}
\sum_{t=1}^{n_0} \log p_{1j}(X_{1j(t-1)},X_{1jt};\theta),~~  
 \wht \theta = {\rm arg}\max_{\theta \in \Theta} L(\theta). 
\end{equation}
Let $k=1$.
%No adjustment of the MLE is performed when the underlying parameter space is finite.

\medskip
\noindent 2. {\it Experimentation}. Let $\lfloor \cdot \rfloor$ denote the greatest integer function.

(a) If $\wht \theta \in \cup_{i > k} \Theta_i$: Take 
$\lfloor z_{kj}(\wht \theta) \log N \rfloor$ observations from job $kj$ for
$j=1,\ldots,J_k$.

\smallskip
(b) If $\wht \theta \in \Theta_k$: Take $\lfloor z_{kj}(\wht \theta)
\log N \rfloor$ observations from job $kj$ for $j \not\in J(\wht 
\theta)$.

\smallskip
(c) If $\wht \theta \in \cup_{i < k} \Theta_i$: Skip experimentation
phase.

\medskip
\noindent 3. {\it Testing}. Start with a full set $\{ k1,\ldots,kJ_k \}$ of
unrejected jobs. Let $\bn = (n_{11},\ldots,n_{kJ_k})$, where $n_{ij}$ denotes
the number of observations taken from arm $ij$ so far.
The rejection of a job is based on the following test statistic. 
Let $F_k$, $1 \leq k \leq I$, be a probability distribution with 
positive probability on all open subsets of $\cup_{i=k}^I \Theta_i$. Define
\begin{equation}\label{teststat}
U_k(\bn;\lambda) = \frac{ \int_{\cup_{i=k}^I \Theta_i} \prod_{i=1}^k
\prod_{j=1}^{J_i} \nu_{ij}(X_{ij0};\theta) \prod_{t=1}^{n_{ij}} p_{ij}(X_{ij(t-1)},X_{ijt};
\theta) \ dF_k(\theta)}{ \prod_{i=1}^k
\prod_{j=1}^{J_i} \nu_{ij}(X_{ij0};\lambda) 
\prod_{t=1}^{n_{ij}} p_{ij}(X_{ij(t-1)},X_{ijt};\lambda) } 
\end{equation}
for all $\lambda \in \Theta_k$. 
 
(a) If $\wht \theta \in \cup_{i > k} \Theta_i$: Add one observation
from each unrejected job. Reject parameter $\lambda$ if $U_k(\bn;\lambda) \geq N$. 
Reject a job $kj$ if all $\lambda \in \Theta_{kj}$ have been rejected at some 
point in the testing stage. If there is a job in group $k$ left unrejected and the total number
of observations is less than $N$, repeat 3(a). Otherwise go to step 4.

\smallskip
(b) If $\wht \theta \in \Theta_k$: Add $n_1$ observations from each
unrejected job $kj$, $j \in J(\wht \theta)$ and one observation from
each unrejected job $kj$, $j \not\in J(\wht \theta)$. Reject a job $kj$ if all 
$\lambda \in \Theta_{kj}$ have been rejected at some point in the testing
phase. If there is a job in group $k$ left 
unrejected and the total number of observations is less than $N$, repeat 3(b). 
Otherwise, go to step 4.

\smallskip
(c) If $\wht \theta \in \cup_{i < k} \Theta_i$: Adopt the procedure of 3(a).

\medskip
\noindent 4. {\it Moving to the next group and termination}. 
The strategy terminates once $N$ observations have been collected. 
Otherwise, if $k < I$, increment $k$ by 1 and go to step 2; if $k=I$,   
select all remaining observations from a job $Ij$ satisfying 
$\mu_{Ij}(\wht \theta) = \max_{1 \leq h \leq J_I} \mu_{Ih}(\wht \theta)$. 

\medskip 
We shall now describe how each feature of the proposed strategy leads to
asymptotic optimality in Theorem 2. The positive information assumption 
in the first half of A2 allows us to estimate $\theta$ consistently and 
hence enables us to determine the optimal sample size $z_{1j}(\theta)$ 
in the experimental stage of group 1. The assumption is important because 
once we move to the next group of jobs, irreversibility would prevent us 
from making up any shortfall in the optimal sample size required 
from group 1. By selecting $n_0 \rightarrow \infty$, we ensure the consistency
of $\wht \theta$ while by choosing $n_0 = o(\log N)$, the estimation of 
$\theta$ incurs negligible contribution to the regret.

Let $k$ be the current group of jobs under sampling.
Consider first $\wht \theta \in \Theta_\ell$ for some $\ell \geq k$. 
We are instructed to select $\lfloor z_{kj}(\wht \theta) \log N \rfloor$
observations from each job in the experimental stage. By Theorem 1 and 
the consistency of $\wht \theta$, this is optimal for learning. 
If $\wht \theta \in \Theta_\ell$ for some $\ell < k$, then the estimate 
$\wht \theta$ says that we have overshot the optimal group, the estimate 
$\wht \theta$ cannot be trusted. In both cases, our strategy then is to rely 
on the testing stage to decide if we should stay within the current job group. 

The testing stage is important in stopping us from moving beyond the
first group of optimal jobs. The rationale is that by irreversibility, 
the penalty for moving beyond the first group of optimal jobs can be of 
order $N$, which is large compared to the desired regret of $O(\log N)$. 
The usefulness of the testing stage in this aspect can be seen from 
(\ref{overshootregret}) 
below, which guarantees that the regret due to overshooting the optimal job
group is $O(1)$. The positive information assumption in the second half of A2 
is necessary for the testing stage to be successful. 
%information can be gathered while sampling from group $k$ and 
%hence allows us to move to the next group if $\theta \in \Theta_\ell$ for some
%$\ell > k$. 

Let us now consider the strategy in 3(b). If $\wht \theta = \theta$, 
then $\lfloor z_{kj}(\theta) \log N \rfloor$ observations from arm $kj$ is
taken in the experimental stage and hence 
by the last inequality of (\ref{mpconstraint}), $o(\log N)$ observations from jobs
with positive information are needed to reject $\lambda \in 
B_\ell(\theta)$ in the testing stage but we may still need an order of $\log N$
observations to reject $\lambda \in \Theta_\ell \setminus B_\ell
(\theta)$. Since we would like $o(\log N)$ observations
from suboptimal jobs in the testing phase, sampling equally from all jobs would
be undesirable here. We consider instead the selection of $n_1$ observations from 
job $\ell j$, $j \in J(\wht \theta)$ for each observation from the
other jobs, where $n_1$ goes to infinity with $N$, so that
$O(n_1^{-1} \log N) = o(\log N)$ observations are taken from suboptimal jobs when 
$\wht \theta = \theta$. When $\wht \theta \neq \theta$, 
it might be possible that each job $k j$, $j \in J(\wht \theta)$
would provide no information to reject some $\lambda \in \Theta_k \setminus 
\cup_{j \in J(\theta)} \Theta_{k j}$. Our procedure
would then allocate $O(n_1 \log N)$ observations from suboptimal jobs in the testing
phase conditional on this happening. By A5 and Chebyshev's inequality, the
probability of providing an incorrect estimate of $\theta$ is
$O(n_0^{-1})$ and hence by specifying $n_1 = o(n_0)$, we 
ensure that the average contribution from suboptimal jobs is  
$O(n_0^{-1} n_1 \log N)=  o(\log N)$.

%Under the finite second moment condition A5 
%and Chebyshev's inequality, the probability of such a scenario occuring is 
%$O(n_0^{-1})$ and hence the overall contribution to the regret for this case is 
%$o(\log N)$ when $n_1=o(n_0)$. 
  
The final case $\wht \theta \in \cup_{i < k} \Theta_i$ occurs with $o(1)$
probability, which together with the $O(\log N)$ observations taken in the
non-optimal jobs in the testing stage when this happens, results in an overall
$o(\log N)$ contribution to the regret. 

The last step is to proceed to the next group of jobs when all
parameters in $\Theta_k$ have been rejected. The exception is when
$k=I$. To be at stage 4 when $k=I$, all $\theta \in \Theta$ have been
rejected at some point in time. Clearly, the true parameter has been
rejected as well but this occurs with very small probability and the 
contribution to the regret in this case is asymptotically negligible. 

\subsection{Extension to infinite $\Theta$}

Let $\theta \in \Theta_\ell$ be the true underlying parameter.
When $\Theta$ is finite, consistency of $\widehat \theta$
would imply that $\widehat \theta = \theta$ with probability close to 1 
when $N$ is large. Hence $B_\ell(\widehat \theta)$ and
$J(\widehat \theta)$ would be good substitutes for the unknown $B_\ell(\theta)$
and $J(\theta)$ respectively. Complications
arise when $\Theta$ is infinite. Firstly, it is possible that
$B_\ell(\theta)$ is non-empty while $B_\ell(\theta')$ is empty for all
$\theta'$ arbitrarily close to $\theta$. Secondly, by 
continuity of $\mu_{ij}(\cdot)$, it follows that there exists
$\delta > 0$ such that 
\[
J(\theta') \subset
J(\theta) {\rm \ for \ all} \ \theta' \in N_\delta(\theta) \cap \Theta_\ell,
\]
but the preceding statement with $\subset$ replaced by $=$ is not
necessarily true. Hence
$B_\ell(\widehat \theta)$ and $J(\widehat \theta)$ are in general poor substitutes of
$B_\ell(\theta)$ and $J(\theta)$ when $\Theta$ is infinite. Moreover if
$\theta$ lies on the boundary of $\Theta_\ell$, then 
$(\cup_{i > \ell} \Theta_i) \cap N_\delta(\theta)$ can be nonempty for all
small $\delta > 0$. This implies that $z_{kj}(\widehat \theta)$ may be inconsistent for
$z_{kj}(\theta)$. This would not happen when $\Theta$ is finite. 

Our strategy in extending the optimal procedure from finite $\Theta$ to infinite
$\Theta$ is not to select $\wht \theta$ during the estimation phase
but rather to select some appropriate adjusted estimate 
$\wht \theta_a \in N_{\delta/2}(\wht \theta)$ where $\delta \rightarrow 0$ 
as $N \rightarrow \infty$ at a rate that is specified in Theorem 2 below. 
We require firstly that
\begin{equation}\label{thetaa1}
\wht \theta_a \in N_{\delta/2}(\wht \theta) \cap \Theta_\ell \quad {\rm where} \ 
\ell = \min \{ i: \Theta_i \cap N_{\delta/2}(\wht \theta) \neq \emptyset \}.
\end{equation}
This condition ensures that if $\theta$ lies in the boundary of $\Theta_\ell$,
then the probability that $\wht \theta_a \in \Theta_\ell$ tends to 1 as $N \rightarrow
\infty$. Our next condition would ensure that the probability that
$J(\wht \theta_a) = J(\theta)$ tends to 1 as $N \rightarrow \infty$.
Let $| \cdot |$ denote the number of elements in a finite set and
$$
{\bf J} = \max \{ |J(\theta')|: \theta' \in N_{\delta/2}(\wht \theta) \cap
\Theta_\ell \}.
$$
We require in addition to (\ref{thetaa1}), that 
\begin{equation}\label{thetaa2}
\wht \theta_a \in H:= \{ \theta \in N_{\delta/2}(\wht \theta) \cap \Theta_\ell:
|J(\theta)| = {\bf J} \}, \mbox{ where } \ell \mbox{ is defined in (\ref{thetaa1}) }.
\end{equation}
If $\Theta$ is finite, then for $\delta > 0$ small
enough, $N_{\delta/2}(\wht \theta) = \{ \wht \theta \}$ and hence by (\ref{thetaa1}) and 
(\ref{thetaa2}), $\wht \theta_a = \wht \theta$. Therefore the selection of $\wht \theta_a$ 
for infinite $\Theta$ is consistent with the selection procedure for finite $\Theta$ 
when $N$ is large. The final thing left to do is the estimation of
$B_\ell(\theta)$. This can be done by taking a union of $B_\ell(\theta')$ over
$\theta' \in H$. We thus have the following modification of the optimal
strategy for infinite $\Theta$, which reduces to the optimal strategy for
finite $\Theta$ for $\delta > 0$ small enough.

\medskip {\bf Optimal strategy for infinite $\Theta$}.

\medskip
1.$'$ {\it Estimation}. Let $k=1$ and $\wht \theta_a$ be an adjusted MLE
satisfying (\ref{thetaa1}) and (\ref{thetaa2}).

\medskip
2.$'$ {\it Experimentation}. 
Let $\wht z_{kj}$ be the solution to Problem A with
parameter $\wht \theta_a$ and with the bad set $B_\ell(\theta)$ replaced
by $\cup_{\theta' \in H} B_\ell(\theta')$. 

(a)$'$ If $\wht \theta_a \in \cup_{i > k} \Theta_i$: Take
$\lfloor \wht z_{kj} \log N \rfloor$ observations from job $kj$,
$j=1,\ldots,J_k$. 

\smallskip
(b)$'$ If $\wht \theta_a \in \Theta_k$: Take $\lfloor \widehat z_{kj}
\log N \rfloor$ observations from job $kj$ for $j \not\in J(\wht \theta_a)$.

\smallskip
(c)$'$ If $\wht \theta_a \in \cup_{i<k} \Theta_i$: Skip experimentation phase.

\smallskip
3$'$. and 4$'$. Identical to the strategy for finite $\Theta$, with $\wht
\theta_a$ replacing $\wht \theta$.

\medskip In view of (\ref{thetaa1}) and (\ref{thetaa2}), 
the modified strategy $\phi^*$ described above will lead
to asymptotic efficiency for infinite $\Theta$ as stated in Theorem 2 below. 
It is also convenient, when $\Theta$ is infinite, to decide on the 
rejection of a job in step 3 based on the current sample rather than 
to keep track of which $\lambda$ have been rejected previously. Hence
for practical use, we can also make the following modification to the
rejection of jobs in step 3$'$:

\smallskip \quad Let $U_{kj}(\bn) = \inf_{\lambda \in \Theta_{kj}} U_k(\bn;\lambda)$.
Reject job $kj$ if $U_{kj}(\bn) \geq N$. 

\begin{theorem} Assume {\rm A1-A5.} 
The strategy $\phi^*$ has error probabilities from the estimation stage
satisfy the following properties.  
Let $n_0 \rightarrow \infty$ with $n_0 = o(\log N)$ and $n_1 \rightarrow \infty$
such that $n_1=o(n_0)$. Then there exists $\delta (=\delta_N) \downarrow 0$ such that
\begin{equation}\label{ep}
P_\theta \{ \wht \theta_a \in \Theta \setminus N_\delta(\theta) \} =o(n_1^{-1}) \ {\rm as} 
\ N \rightarrow \infty.
\end{equation}
Let $\theta \in \Theta_\ell$. Then the regret of $\phi^*$ due to overshoot in the testing stage 
is $O(1)$ because
\begin{equation}\label{overshootregret}
\sum_{i>\ell} \sum_{j=1}^{J_i} E_{\theta} T_N(ij) \leq 1.
\end{equation}
Therefore, the total regret   
\begin{eqnarray}\label{optimalregret}
\lim_{N \rightarrow \infty} R_{N}(\theta)/\log N = z(\theta,\ell).
\end{eqnarray}
%In addition, the average sampling cost
%\begin{equation}
%H_N(\theta) = o(\log N) \mbox{ as } N \rightarrow \infty.
%\end{equation}
\end{theorem}
\noindent{\bf Remark 1.} Theorem 2 extends Fuh and Hu (2000) to situations where 
more than one job in each group are available for processing. Theorem 2 generalizes 
the results of Lai (1987) and Agrawal et al. (1989a,b) to the case of infinite state
and parameter spaces and more than one job group. 

\noindent{\bf Remark 2.} If there is a constant switching cost each time we switch 
from one job to another, it can be shown that the strategy $\phi^*$ has switching 
cost of order $o(\log N)$. Hence
$\phi^*$ is still efficient considering switching cost. The details will be given in another paper. 

\noindent{\bf Remark 3.} We consider non-empty bad set in this paper. It can be shown that the proposed
strategy $\phi^*$ can achieve $o(\log N)$ regret, when the bad set is empty and $I=1$. 
In general, within the optimal group, the contribution to the regret from jobs optimal 
for parameter values outside the
bad set is $o(\log N)$. The essence of the proof for this fact is contained in 
Section 6. We will  provide detailed justification in another paper. The upshot is 
that the strategy $\phi^*$ can achieve
super efficient results outside of bad sets.

%3. In Agrawal, Hedge and Teneketzis (1988) and Graves and Lai (1997),
%which considers one group of correlated arms, (4.12) is also achieved
%by sampling in blocks. 

\section{Examples}
\def\theequation{5.\arabic{equation}}
\setcounter{equation}{0}

\noindent {\bf Example 3: Multi-phase project management.}
To illustrate how our method can be applied, we discuss a few 
examples. Our purpose here is not to provide an accurate statistical 
model for a particular situation, but rather to supply concrete examples
of parameter spaces and probability distributions such that the assumptions in 
Section 2.3 are satisfied. 

Consider the management of $N$ research and development (R\&D) projects. 
When a project is pursued with one unit of resource, the reward is a normal 
random variable $X$ with mean $\mu_t(\theta)$ and variance $\sigma^2_t(\theta)$. 
Given the parameter value $\theta$, the mean 
$\mu_t(\theta)$ reflects, at time $t$, the level of existing technology and 
knowledge relevant to the concerned projects as well as the competition in the 
market. Let $\theta=(\alpha,\beta)$ and 
\begin{equation}\label{msd}
\mu_t(\theta) = {f(t,\alpha) \over h(t,\beta)}, \qquad \sigma_t(\theta) = 
{1 \over h(t,\beta)}, 
\end{equation}
where both $f(t,\alpha)$ (reflecting technology and knowledge) and $h(t,\beta)$ 
(reflecting competition) are increasing functions of time $t$. Observe that under 
(\ref{msd}) the coefficient of variation, 
$\sigma_t/\mu_t=1/f(t,\alpha)$ is a decreasing function of $t$, which can be 
interpreted as follows. Because the products from the project will be gradually 
superseded by more advanced ones through competition in the market, therefore not 
only the mean reward becomes smaller but we are also more certain of it due to 
as time moves on. If we take $f$ and $h$ to be
\begin{equation}\label{examplemsd}
f(t,\alpha) = \alpha t^2, \quad {\rm and} \quad h(t,\beta) = e^{t \beta}-1,
\end{equation}
then the maximal value of $\mu_t(\theta)$, for a fixed value of $\theta$, 
is attained uniquely at $t$ such that $t\beta={\rm constant}\approx 1.5936$. 

Designate $I$ phases indexed by time points $0 < t_1 < \cdots < t_I$ during which 
pursuing a project can take place. And there are $J$ different types of projects that can be pursued 
at any phase $i=1,\ldots,I$. To accommodate $I$ phases and $J$ types of projects, we expand the 
parameter vector to $\theta = (\alpha_1,\ldots,\alpha_J,\beta)$. 
Given (\ref{msd}) and (\ref{examplemsd}), let the reward $X_{ijk}$ from the pursue with $k$-th 
unit of resource of the type $j$ project in phase $i$ be i.i.d. normal with means and standard deviations
\begin{equation}\label{msdij}
\mu_{ij}(\theta) = {\alpha_j t_i^2 \over e^{t_i \beta}-1}, \qquad \sigma_{i}(\theta) = 
{1 \over e^{t_i \beta}-1}, 
\end{equation}
respectively. 

By selecting
$\Theta = [\underline{\alpha},\overline{\alpha}]^J \times [\underline{\beta},
\overline{\beta}]$ where $0 < \underline{\alpha} < \overline{\alpha} < \infty$
and $0 < \underline{\beta} < \overline{\beta} < \infty$, condition A1 is
easily seen to hold. Let $\theta'=(\alpha'_1,\cdots,\alpha'_J,\beta')$, then 
$$I_{ij}(\theta,\theta')=\log \Big[ \frac{\sigma_i(\theta')}{\sigma_i(\theta)} \Big]+
\frac{\sigma_i^2(\theta)-\sigma_i^2(\theta')+[\mu_{ij}(\theta)-\mu_{ij}(\theta')]^2}{2\sigma_i^2(\theta')} 
$$
equals zero if and only if
$\mu_{ij}(\theta) = \mu_{ij}(\theta')$ and $\sigma_{ij}^2(\theta) = 
\sigma_{ij}^2(\theta')$, or equivalently, $\alpha_j  =\alpha_j'$ and
$\beta = \beta'$, the information assumption A2 is also satisfied. 
It can be shown that there exist $(\underline{\beta}=) \beta_I < \beta_{I-1} < \cdots < 
\beta_1 < \beta_0 (=\overline{\beta})$ such that
\begin{equation}\label{Thetai}
\Theta_i = \cases {\{ \theta \in \Theta: \beta \in [\beta_1,
\beta_0] \} \ {\rm for} \ i=1 \cr
\{ \theta \in \Theta: \beta \in [\beta_i,\beta_{i-1}) \}
\ {\rm for} \ 2 \leq i \leq I, }
\end{equation}
and 
\begin{equation}\label{Thetaij}
\Theta_{ij} = \{ \theta \in \Theta_i : \alpha_j = \max_{1 \leq k
\leq J} \alpha_k \}.
\end{equation}

Since the observations $X_{ijk}$ are independent, the assumptions
A3-A5 are satisfied by selecting $G_{ij} = {\bf R}$, $V_{ij}(x)=|x|+1$ and
$b > \sup_{i,j,\theta} E_\theta |X_{ij1}|+1$.
Consequently, the strategies described in Section 4 
are efficient in the sense of attaining the regret lower bound given by Theorem 1.\\ 

\noindent {\bf Example 4: Multi-phase project management with Markovian reward.} 
Continuing from Example 3, instead of i.i.d. reward, we assume that $k$-th pursue 
of a project of type $j$ at time $t_i$ follows an AR(1) process 
\[X_{ijk} = a_i X_{ij(k-1)} + \epsilon_{ijk},\]
where $|a_i|<1$ and $\epsilon_{ijk}\sim N(\mu_{ij}(\theta), \sigma_{i}^2(\theta))$ with 
$\mu_{ij}$ and $\sigma_i$ given by (\ref{msdij}). Let $G_{ij}=[-c,c]$ for some
$c > 0$. Since $\epsilon_{ij1}$ has a positive density on
the real line, A3 is satisfied. Let $V_{ij}(x)=|x|+1$. 
From Meyn and Tweedie (1993) page 380, $\{X_{ijk}\}_{k \geq 0}$ is geometric 
ergodic and A4 holds with $0 < \bar{b} < 1-\max_i |a_i|$ and $b,c$ large enough.

The stationary distribution is normally distributed with mean and variance given by 
$(1-a_i)^{-1}\mu_{ij}(\theta)$ and $(1-a_i^2)^{-1}\sigma_i^2(\theta)$.  It can be 
checked that (\ref{Thetai}) and (\ref{Thetaij}), which reveal the structure of the 
parameter space, still holds for AR(1) reward. Consequently, 
A1 is true for AR(1) rewards. To simplify the presentation of the Kullback-Leibler information number,
we drop the indices $i,j$ and use $\mu',\sigma'$ to denote $\mu(\theta'),\sigma(\theta')$, respectively.  
\begin{eqnarray*}
I(\theta,\theta')&=&
\log \Big( \frac{\sigma'}{\sigma} \Big) +
\frac{\sigma^2-\sigma'^2+(\mu-\mu')^2}{2\sigma'^2}+\\
&+&\frac{(a-a')^2\{\mu^2(1-a)^{-2}+\sigma^2(1-a^2)^{-1}\}+2(a-a')(\mu-\mu')\mu(1-a)^{-1}}{2\sigma'^2}.
\end{eqnarray*}
It is clear that the Kullback-Leibler number is greater than zero if $\theta\neq\theta'$. From the preceding 
equation, we can verify that A2 and A5 hold.

%\begin{eqnarray*}
%I_{ij}(\theta,\theta')&=&
%\log\frac{\sigma_i(\theta')}{\sigma_i(\theta)}+
%\frac{\sigma^2_i(\theta)-\sigma^2_i(\theta')+[\mu_{ij}(\theta)-\mu_{ij}(\theta')]^2}
%{2\sigma^2_i(\theta')}+\\
%&+&\frac{(\alpha_j-\alpha'_j)^2\{\mu^2(1-\alpha_j)^{-2}+\sigma_i^2(1-\alpha_j^2)^{-1}\}+
%2(\alpha_j-\alpha'_j)(\mu_{ij}(\theta)-\mu_{ij}(\theta)')\mu_{ij}(1-\alpha_j)^{-1}}
%{2\sigma_i^2(\theta')},
%\end{eqnarray*}

%\subsection{Computerize adaptive tests}

\section{Proof of asymptotic efficiency}
\def\theequation{6.\arabic{equation}}
\setcounter{equation}{0}

We shall demonstrate the asymptotic efficiency of $\phi^*$ by proving
(\ref{ep})-(\ref{optimalregret}). A change-of-measure argument is first used to prove 
(\ref{overshootregret}). As the proofs of (\ref{ep}) and (\ref{optimalregret}) are too involved for one
reading, we prove them in Section 6.1 for the restricted case of finite $\Theta$
and extend the proofs to infinite $\Theta$ in Section 6.2.

\medskip \noindent {\bf Proof of (\ref{overshootregret})}. 
Let $\wtd P$ be the measure which generates ${\bf X}_\bn := \{ X_{ijt} \}$ for
$j=1,\ldots,J_i$ and $i=1,\ldots,\ell$, $t=1,\ldots,n_{ij}$
in the following manner. First generate $\theta'$ randomly from
$F_{\ell}$. Using the strategy $\phi^*$ to select the jobs to be processed, 
generate $X_{ij0}$ from $\nu_{ij}(\cdot;\theta')$ and $X_{ijt}$, $t \geq 1$, according
to the transition density $p_{ij}(X_{ij(t-1)},\cdot;\theta')$ when 
at job $ij$. Let $\theta \in \Theta_{\ell j}$. Then
\[
\frac{d \wtd P}{dP_{\theta}}({\bf X}_\bn) = 
\frac{ \int_{\cup_{i=\ell}^I \Theta_i} \prod_{i=1}^\ell
\prod_{j=1}^{J_i} \nu_{ij}(X_{ij0};\theta') \prod_{t=1}^{n_{ij}} p_{ij}(X_{ij(t-1)},X_{ijt};
\theta') \ F_\ell(d \theta')}{ \prod_{i=1}^\ell
\prod_{j=1}^{J_i} \nu_{ij}(X_{ij0};\theta) 
\prod_{t=1}^{n_{ij}} p_{ij}(X_{ij(t-1)},X_{ijt};\theta)} = U_{\ell}(\bn;\theta).
\]
Let $\bT = (T_N(11),\ldots,T_N(\ell J_{\ell}))$ and $A = 
\{ U_{\ell}(\bT;\theta) \geq N \}$. Then $P_\theta \{ \sum_{i > \ell} T_N(i) > 0 \}$
is bounded by 
\begin{equation}\label{overshootprobability}
P_{\theta}(A) = E_{\tilde{P}} \Big[\frac{dP_{\theta}}{d \wtd P} ({\bf X}_{\bT}) {\bf 1}_A \Big] \leq
N^{-1}. 
\end{equation}
Hence (\ref{overshootregret}) follows from (\ref{overshootprobability}) and the bound 
$\sum_{i>\ell} T_N(i) \leq N$. $\Box$
 
\subsection{Finite parameter space}

Let $\Theta = \{ \theta_0,\ldots,\theta_h \}$.  Let $\theta_0 \in \Theta_{\ell j_0}$ 
be the true  parameter value. For $1 \leq q \leq h$, define 
\begin{equation}\label{xiq}
\xi_{ijt}(q) = \log [p_{ij}(X_{ij(t-1)},X_{ijt};\theta_0)/
p_{ij}(X_{ij(t-1)},X_{ijt};\theta_q)].
\end{equation}
Then $E_{\pi(\theta_0)} \xi_{ijt}(q) = I_{ij}(\theta_0,\theta_q)$.  To get the 
essence of the strategy without being overly involved in cumbersome notation, let us consider a specific case 
$\ell=2$, $J_1=J_2=2$, $\theta_0 \in \Theta_{21}$ and $J(\theta_0) = \{ 1 \}$. 

We first prove (\ref{ep}). Let us consider the inequality
\begin{eqnarray*}
& & P_{\theta_0} \{ \wht \theta \neq \theta_0 \} = \sum_{q=1}^h P_{\theta_0}
\{ \wht \theta = \theta_q \} 
\leq \sum_{q=1}^h  P_{\theta_0} \Big\{ \sum_{t=1}^{n_0}
\xi_{11t}(q) + \xi_{12t}(q) < 0 \Big\} . 
\end{eqnarray*}
By A5 and Chebyshev's inequality, 
\begin{eqnarray*}
& & P_{\theta_0} \Big\{ \sum_{t=1}^{n_0} \xi_{11t}(q) + \xi_{12t}(q)< 0 \Big\}
\leq { {\rm Var}_{\theta_0} \Big( \sum_{t=1}^{n_0} \xi_{11t}(q)+ \xi_{12t}(q)\Big)
\over \Big[ E_{\theta_0} \Big( \sum_{t=1}^{n_0} \xi_{1jt}(q)+\xi_{11t}(q) \Big)
\Big]^2} \\
& & \quad \leq (1+o(1)) {E_{\pi(\theta_0)} [\xi_{11t}^2(q)+\xi_{12t}^2(q)]+
2I_{11}(\theta_0,\theta_q)I_{12}(\theta_0,\theta_q) \over
n_0 [I_{11}(\theta_0,\theta_q)+I_{12}(\theta_0,\theta_q)]}=O(n_0^{-1}).
\end{eqnarray*}
This completes the proof of (\ref{ep}) for finite parameter case. 

We now undertake the proof of (\ref{optimalregret}). For $q \geq 1$, let
$\theta_q \in \Theta_{k j'}$ where either (i) $k < \ell$ or (ii) $k=\ell$
and $j' \notin J(\theta_0)$. Let $\tau_{kj}(q)$ be the number of observations 
selected from job $kj$ in the testing phase of group $k$ before parameter 
$\theta_q$ is rejected. To show (\ref{optimalregret}), it suffices to prove that
\begin{equation}\label{testingregret}
E_{\theta_0} \Big[ \sum_{j=1}^{J_k} \tau_{kj}(q) \Big] = o(\log N) \ {\rm if} \
k < \ell ~~{\rm and}~~ \ E_{\theta_0} \Big[ \sum_{j \notin J(\theta_0)} 
\tau_{\ell j}(q) \Big]  = o(\log N)
\end{equation}
because (\ref{testingregret}) implies that the regret in the testing phase before 
leaving the optimal group is $o(\log N)$ and the regret due to overshooting the optimal group, which is also 
$o(\log N)$ by the established (\ref{overshootregret}), complete the justification. 

%hence $H_N(\theta_0)=o(\log N)$ because there is $O(1)$ switching 
%of arms in both the estimation and experimentation phase.

Select $C > 0$ large enough such that $\xi_{ijt}'(q) = \xi_{ijt}(q) \wedge C$
has positive expectation under $\pi(\theta_0)$ for all $i,j,q$ 
satisfying $I_{ij}(\theta_0,\theta_q)$ $> 0$. Let $\bn = (n_{11},n_{12})$. We will first show that the first half of (\ref{testingregret}) is satisfied when $\theta_q \in \Theta_1$. By (\ref{teststat}),
\begin{equation}\label{logteststat}
\log U_1(\bn;\theta_q) \geq \sum_{t=1}^{n_{11}} \xi_{11t}(q) +
\sum_{t=1}^{n_{12}} \xi_{12t}(q) + \log v_{11} + \log v_{12} + \log F_1(\theta_0),
\end{equation}
where $v_{ij} = \inf_{x,\theta,\lambda} [\nu_{ij}(x;\theta)/\nu_{ij}(x;\lambda)] > 0$
as assumed in (\ref{2.2a}). Hence by (\ref{teststat}), rejection of $\theta_q$ has 
occurred when
\begin{eqnarray} \label{definec}
& \sum_{t=1}^{m_{11}} \xi_{11t}(q) + \sum_{t=1}^{m_{12}} \xi_{12t}(q) 
+ \sum_{t=m_{11}+1}^{n_{11}} \xi_{11t}'(q) + \sum_{t=m_{12}+1}^{n_{12}}
\xi_{12t}'(q) \cr
& \quad > c := \log N - \log
v_{11} -\log v_{12} - \log F_1(\theta_0),
\end{eqnarray}
where $\bm = (m_{11},m_{12})=(n_0+\lfloor z_{11}(\wht \theta) \log N \rfloor,
n_0+\lfloor z_{12}(\wht \theta) \log N \rfloor)$ is the
sample size at the beginning of the testing phase. Since
$\xi_{ijt}'(q)$ is bounded above by $C$, it follows that at $\bn=(n'_{11},n'_{12})$ 
for which the boundary is first crossed by $\xi'_{ijt}$'s
\begin{equation}\label{crossing} \qquad
E_{\theta_0} \Big[ \sum_{t=1}^{m_{11}} \xi_{11t}(q) + \sum_{t=1}^{m_{12}}
\xi_{12t}(q) \Big] + 
E_{\theta_0} \Big[ \sum_{t=m_{11}+1}^{n'_{11}} \xi_{11t}'(q) + 
\sum_{t=m_{12}+1}^{n'_{12}}
\xi_{12t}'(q) \Big] \leq c(1+o(1)).
\end{equation}

By (\ref{ep}), the condition $n_0 = o(\log N)$, (\ref{definec}), and the constraint 
$I_{11}(\theta_0,\theta_q) z_{11}(\theta_0) + I_{12}(\theta_0,
\theta_q)$ \ $z_{12}(\theta_0) \geq 1$ from (\ref{mpconstraint}), it follows that
\begin{equation}\label{tsss}
E_{\theta_0} \Big[ \sum_{t=1}^{m_{11}} \xi_{11t}(q) + \sum_{t=1}^{m_{12}}
\xi_{12t}(q) \Big] \geq (1+o(1)) c.
\end{equation}
Subtracting (\ref{tsss}) from (\ref{crossing}), we have
\begin{equation}\label{tsincrement}
E_{\theta_0} \Big[ \sum_{t=m_{11}+1}^{n_{11}} \xi_{11t}'(q) + 
\sum_{t=m_{12}+1}^{n_{12}}
\xi_{12t}'(q) \Big] = o(c).
\end{equation}
By Wald's equation for Markov processes, the left hand side of 
(\ref{tsincrement}) equals
\begin{equation}\label{tsincrement'}
(1+o(1)) [E_{\pi(\theta_0)} \xi_{11t}'(q) E_{\theta_0}(n_{11}-m_{11}) +
E_{\pi(\theta_0)} \xi_{12t}'(q) E_{\theta_0} (n_{12}-m_{12})].
\end{equation}
%For independent rewards, we have strict inequality instead of a $(1+o(1))$equivalence. 
The proof of Wald's equation for Markov process is given in
the Appendix. By A2 and the choice of $C$ for $\xi'_{ijt}(q)$, 
$E_{\pi(\theta_0)} \xi_{1jt}'(q) > 0$.  In view of the sample size in testing stage 
$\tau_{1j}(q) \leq n_{1j} -m_{1j}$, it follows from 
(\ref{tsincrement})-(\ref{tsincrement'}) that 
$E_{\theta_0} [\tau_{11}(q)+\tau_{12}(q)] = o(c)$
for all $\theta_q \in \Theta_1$. Hence the rejection of both
$\Theta_{11}$ and $\Theta_{12}$ involves only $o(\log N)$ observations and the first half of (\ref{testingregret}) holds. 

Next we show that the second half of (\ref{testingregret}) holds 
when $\theta_q \in \Theta_{22}$. We divide into two cases, $\theta_q \in B_2(\theta_0)$ and
$\theta_q \not\in B_2(\theta_0)$. Consider the first case. By (\ref{bsdual}),
$I_{22}(\theta_0,\theta_q) > 0$. We then follow the arguments above using (\ref{ep}) and the last inequality of (\ref{mpconstraint}) to show that $E_{\theta_0} \tau_{22}(q) = o(c)$. 

The second scenario involves $\theta_q \not\in B_2(\theta_0)$.
The key observation is $I_{21}(\theta_0,\theta_q) > 0$ by (\ref{bs}). 
In other words, information is always collected and no additional
regret is incurred when we sample from job 21. Under unequal sampling,
\begin{eqnarray}\label{tau22}
E_{\theta_0} \tau_{22}(q) & = & E_{\theta_0}[\tau_{22}(q) {\bf 1}_{\{
J(\hat \theta)=\{1 \} \}}] + E_{\theta_0}[\tau_{22}(q) {\bf 1}_{\{
J(\hat \theta)=\{1,2 \} \}}] + E_{\theta_0}[\tau_{22}(q) {\bf 1}_{\{
J(\hat \theta)=\{2 \} \}}] \cr
& = & n_1^{-1} E_{\theta_0} [ \tau_{21}(q) {\bf 1}_{\{ J(\hat \theta)= 
\{1 \} \}}]
+ E_{\theta_0} [\tau_{21}(q) {\bf 1}_{\{ J(\hat \theta) = \{1,2 \} \}}] \cr
& & \qquad + n_1 E_{\theta_0} [\tau_{21}(q) {\bf 1}_{\{ J(\hat \theta)
=\{ 2\} \}}]. 
\end{eqnarray}
Since $
E_{\theta_0} \Big[ \tau_{21}(q) {\bf 1}_{\{ J(\hat 
\theta) = A \}} \Big] \leq (1+o(1))c P_{\theta_0} \{ J(\hat \theta) = A \}/
I_{21}(\theta_0,\theta_q)$
for $A = \{1 \}, \{ 2 \}$ and $\{ 1,2 \}$ and as $n_1 \rightarrow\infty$, the first term on 
the right hand side of (\ref{tau22}) is $o(c)$ while (\ref{ep}) ensures the second term is $o(c)$. 
By (\ref{ep}), $n_1 P_{\theta_0} \{ J(\wht \theta) \neq \{1 \} \} \leq n_1 P_{\theta_0}
\{ \wht \theta \neq \theta_0 \} = o(1)$
and thus the third term on the right hand side of (\ref{tau22}) 
is $o(c)$. We can conclude that $E_{\theta_0} \tau_{22}(q)=o(c)$ or the second half of 
(\ref{testingregret}) holds. 

\subsection{Extension to infinite parameter space}
We preface the extension with the following lemma. 
The proof of this lemma is given in the Appendix in Section 7.3.
We shall let $\bar A$ denote the closure of a set $A$.

\begin{lemma} Let $\theta_0 \in \Theta_{\ell}$. Assume {\rm A1-A5} and let
$n_0 \rightarrow \infty$, $n_1=o(n_0)$. 

\smallskip 
\noi {\rm (a)} \quad Let $\theta' \not= \theta_0$ and let 
$\wht \theta$ be the MLE estimate {\rm (\ref{mle})}. Then there exists
$\delta' > 0$ small enough such that 
\begin{equation}\label{epmle}
P_{\theta_0} \{ \widehat \theta \in N_{\delta'}(\theta') \} \rightarrow
0 \ as \ N \rightarrow \infty.
\end{equation} 
{\rm (b)} \quad Let $\theta' \in \bar \Theta_k$ for some $k < \ell$
or $\theta' \in \cup_{j \not\in J(\theta_0)} \bar \Theta_{\ell j}$. 
Let $\delta' > 0$ and let $\tau_{kj}$ ($\tau_{\ell j}$)
be the number of observations selected from job $kj$ ($\ell j$) in the
testing phase of group $k$ ($\ell$) before all parameters in the set
$N_{\delta'}(\theta')$ are rejected. Then for $\delta' > 0$ small
enough, 
\begin{equation}\label{testingregretinf}
E_{\theta_0} \Big( \sum_{j=1}^{J_k} \tau_{kj} \Big)
= o(\log N) \mbox{ if } k < \ell \ \mbox{ and } E_{\theta_0} \Big(
\sum_{j \not\in J(\theta_0)} \tau_{\ell j} \Big) = o(\log N).
\end{equation}
\end{lemma} 

We now apply Lemma 2 to extend the proof of Theorem 2. 
By the compactness of $\Theta \setminus N_{\delta/2}(\theta_0)$, $\delta>0$, there exists a
finite set $\{ \theta_1,\ldots,\theta_h \}$ and constants $\delta_q > 0$
such that (\ref{epmle}) holds for $\theta' = \theta_q$ and $\delta' = \delta_q$
for all $1 \leq q \leq h$ and 
$\Theta \setminus \{ \theta_0 \} \supset 
\bigcup_{q=1}^h N_{\delta_q}(\theta_q) \supset \Theta \setminus N_{\delta/2}(\theta_0)$.
Then by (\ref{epmle}), $P_{\theta_0} \{ \widehat \theta \in \Theta \setminus 
N_{\delta/2}(\theta_0) \} \rightarrow 0$
as $N \rightarrow \infty$ and the result (\ref{ep}) follows from
(\ref{thetaa1}) because $\| \wht \theta_a - \wht \theta \| < \delta/2$. 

It remains to show that the number of observations taken from each non-optimal job 
in the testing phase is $o(\log N)$. Consider $k < \ell, j=1,\ldots,J_k$ or $k=\ell$ with
$j \not\in J(\theta_0)$. Since $\bar \Theta_{k j}$ is compact, 
there exists a finite set
$\{ \theta_1,\ldots,\theta_h \}$ and constants $\delta(q) > 0$ such that (\ref{testingregretinf})
is satisfied for $\theta'=\theta_q$, $\delta'=\delta_q$ for all $1 \leq q \leq h$ and 
$\bigcup_{q=1}^h N_{\delta_q}(\theta_q) \supset \bar \Theta_{kj}$, 
and hence by (\ref{testingregretinf}), the number of times job $kj$ is processed
in the testing phase is $o(\log N)$ as required. 

\def\theequation{7.\arabic{equation}}
\setcounter{equation}{0}
\section{Appendix}
\subsection{Proof of (\ref{2.3a})}

Let $X_{ijt}$ denotes the $t$th observation taken from arm $ij$. Then
\begin{equation} \label{7.1}
\Big| W_N(\theta) - \sum_{i=1}^I \sum_{j=1}^{J_i} \mu_{ij}(\theta) E_\theta
T_N(ij) \Big| \leq \sum_{i=1}^I \sum_{j=1}^{J_i} \sum_{t=1}^\infty
| E_\theta g(X_{ijt}) - \mu_{ij}(\theta) |. 
\end{equation}
For any signed measure $\lambda$ on $(D,{\cal D})$, let
\begin{equation} \label{7.2}
\| \lambda \|_{V_{ij}} = \sup_{h:|h| \leq V_{ij}} \Big| \int h(x) \lambda(dx) \Big|.
\end{equation}
It follows from Meyn and Tweedie (1993, p.367 and Theorem 16.0.1) that under
A3 and the geometric drift condition (\ref{drift}), 
\begin{equation} \label{7.3}
\omega_{ij} := \sup_{\theta \in \Theta, x \in D} \sum_{t=1}^\infty
\| P_{ijt}^\theta(x,\cdot) - \pi_{ij}(\cdot;\theta) \|_{V_{ij}}/ V_{ij}(x) < \infty,
\end{equation}
where $P_{ijt}^\theta(x,\cdot)$ denotes the distribution of $X_{ijt}$ conditioned
on $X_{ij0}=x$ and $\pi_{ij}(\cdot; \theta)$ denotes the stationary distribution of
$X_{ijt}$ under parameter $\theta$. By (\ref{2.10a}), there exists $\kappa > 0$
such that $\kappa |g(x)| \leq V_{ij}(x)$ for all $x \in D$ and hence it follows
from (\ref{7.2}) and (\ref{7.3}) that
\begin{equation} \label{7.4}
\kappa \sum_{t=1}^\infty |E_{\theta,x} g(X_{ijt})-\mu_{ij}(\theta)| \leq
\omega_{ij} V_{ij}(x), 
\end{equation}
where $E_{\theta,x}$ denotes expectation with respect to $P_\theta$ and
intial distribution $X_{ij0}=x$. 

In general, for any initial distribution $\nu_{ij}(\cdot;\theta)$,
it follows from (\ref{vij}) and (\ref{7.4})
that
$$\sum_{t=1}^\infty |E_\theta g(X_{ijt}) - \mu_{ij}(\theta)| \leq
\int \sum_{t=1}^\infty | E_{\theta,x} g(X_{ijt}) - \mu_{ij}(\theta)| \nu_{ij}
(x;\theta) Q(dx) < \infty
$$
uniformly over $\theta \in \Theta$ and hence (\ref{2.3a}) follows from (\ref{7.1}).

\subsection{Proof of Lemma 1}

To prove (\ref{ml}), it suffices to show that for every
$\theta' \in \Theta_k^*$, $k < \ell$ and for $\delta > \alpha >0$,
\begin{equation}\label{mlprob}
\lim_{N\rightarrow \infty} P_{\theta} \Big\{ \sum_{i=1}^k \sum_{j=1}^{J_i}
I_{ij} (\theta,\theta') T_N(ij) < (1-\delta) \log N \Big\} = 0.
\end{equation}
Because $\phi$ is uniformly good and $\theta'\in \Theta_k^*$, it follows from
(\ref{uniformlygood}) that $E_{\theta^{'}}
[N-\sum_{j \in J(\theta')} T_N(kj)] = o(N^{\alpha })$ for $\alpha > 0$.
By A2, $I_{kj}(\theta,\theta') > 0$ for all $j \in J(\theta')$ and hence
$I_0:= \min_{j \in J(\theta')} I_{kj}(\theta,\theta') > 0$. 
It then follows from Chebyshev's inequality that 
\begin{eqnarray}\label{mlprobtheta'}
&~& P_{\theta'} \Big\{ \sum_{i=1}^k \sum_{j=1}^{J_i} I_{ij}
(\theta,\theta') T_N(ij) < (1-\delta) \log N \Big\} \\
& \leq & P_{\theta'} \Big\{ I_0 \sum_{j \in J(\theta')} T_N(kj) <
(1-\delta) \log N \Big\} \nonumber \cr
& = & P_{\theta'} \Big\{ [N -\sum_{j \in J(\theta')} T_N(kj)] > 
N-(1-\delta) (\log N)/I_0 \Big\} \nonumber \cr
&=& O(N^{-1}) E_{\theta^{'}}\Big[ N- \sum_{j \in J(\theta')}
T_N(kj) \Big] = o(N^{\alpha -1}). \nonumber 
\end{eqnarray}

Let $\bn = (n_{11},\ldots,n_{kJ_k})$ and $\bT_N = (T_N(11),\ldots,T_N(kJ_k))$.
Let
$$
L_{\bn}=\sum_{i=1}^k \sum_{j=1}^{J_i} \Big\{ \log [ \nu_{ij}(X_{ij0};\theta)/ 
\nu_{ij}(X_{ij0};\theta') ] + \sum_{t=1}^{n_{ij}} \ell_{ij}
(X_{ij(t-1)},X_{ijt};\theta,\theta') \Big\}
$$
be the log likelihood ratio of $\theta$ with respect to $\theta'$, and denote
\[
G_N=\Big\{\sum_{i=1}^k \sum_{j=1}^{J_i} I_{ij}(\theta,\theta')T_N(ij)
< (1-\delta) \log N~{\rm and}~L_{\bT} \leq
(1-\alpha) \log N \bigg \}.
\]
Then by (\ref{mlprobtheta'}), $P_{\theta^{'}}(G_N)=o(N^{\alpha-1})$. By Wald's likelihood ratio 
identity for Markov chains, 
\begin{eqnarray*}
&~&  P_{\theta'}\left\{\bT_N = \bn, L_{\bn} \leq (1-\alpha) \log N \right\}
= E_\theta \Big[ \exp(-L_{\bn}) {\bf 1}_{\{ \bT_N=\bn, L_{\bn} \leq (1-\alpha) \log N \}}
\Big] \cr
& & \qquad  \geq N^{\alpha-1} P_{\theta} \{ \bT_N=\bn,
L_{\bn} \leq (1-\alpha) \log N\}.
\end{eqnarray*}
By summing the preceding inequality over all $\bn$, we have
\begin{equation}\label{mllrprob}
P_{\theta}(G_N)\leq N^{1-\alpha} P_{\theta^{'}} (G_N)=N^{1-\alpha}
o(N^{\alpha-1})=o(1).
\end{equation}

By A3 and the strong law of large numbers for Markov chains (cf. Theorem 17.0.1 of
Meyn and Tweedie, 1993), 
\begin{eqnarray*}
 \Big|L_{\bf n} - \sum_{i=1}^k \sum_{j=1}^{J_i}
 I_{ij}(\theta,\theta') n_{ij} \Big| =o \Big(\sum_{i=1}^k
 \sum_{j=1}^{J_i} n_{ij} \Big)~~~P_{\theta}~\hbox{a.s. as }
\sum_{i=1}^k \sum_{j=1}^{J_i}
n_{ij} \rightarrow \infty.
\end{eqnarray*}
Thus, 
\begin{eqnarray*}
\lim_{m \rightarrow \infty} \Big\{ 
\max_{\bn: \sum_{i=1}^k \sum_{j=1}^{J_i} I_{ij}(\theta,\theta') n_{ij}
\leq m} \Big[ L_{\bn}-\sum_{i=1}^k \sum_{j=1}^{J_i}
I_{ij}(\theta,\theta')n_{ij} \Big] \Big/ m \Big\} 
\rightarrow 0 \mbox{ a.s. under } P_\theta.
\end{eqnarray*}
Because $1-\alpha>1-\delta$, it then follows that as $N \rightarrow \infty$,
$$
P_{\theta} \Big\{L_{\bn} >(1-\alpha) \log N,~\hbox
{for~some}~\bn~\hbox{such~that} \ \sum_{i=1}^k \sum_{j=1}^{J_i} I_{ij}(\theta,\theta')
    n_{ij} < (1-\delta) \log N \Big\}  \rightarrow 0.
$$
Therefore, as $N \rightarrow \infty$,
\begin{eqnarray*}
 P_{\theta}\Big\{\sum_{i=1}^k \sum_{j=1}^{J_i} I_{ij}(\theta,\theta')
 T_N(ij) < (1-\delta) \log N~~\hbox{and}~~L_{\bT}
 >(1-\alpha) \log N \Big\} \rightarrow  0.
\end{eqnarray*}
This combined with (\ref{mllrprob}) gives (\ref{mlprob}), from which (\ref{ml}) follows by
letting $\delta \downarrow 0$.

We now consider the case $\theta' \in B_{\ell}(\theta)$. By
(\ref{bsdual}), $\min_{j \in J(\theta')} I_{\ell j}(\theta,\theta')>0$.
The proof proceeds as before with $k=\ell$, which leads us to (\ref{ml})
with $k=\ell$. Since $I_{\ell j}(\theta,\theta') = 0$ for all $j \in J(\theta)$
by (\ref{bs}),  (\ref{mlbs}) follows.

\subsection{Proof of Theorem 1}

As we mentioned after (\ref{bsregret}) that $B_{\ell}(\theta)\neq \emptyset$, by A1,
$\Lambda_\ell = \Theta_1^* \times \cdots \times \Theta_{\ell-1}^* 
\times B_\ell(\theta)$ is non-empty. For each $\lambda = (\lambda_1,\cdots,\lambda_\ell)\in 
\Lambda_\ell$ and $\theta \in \Theta_\ell$, we define $z(\theta,\ell,\lambda)$
to be the minimal value of (\ref{mp}) with (\ref{mpconstraint}) replaced by
\begin{equation}\label{mplemma1}
\left\{ \begin{array}{ll}
\sum_{j=1}^{J_1} I_{1j}(\theta,\lambda_1) z_{1j}(\theta) \ge 1,    \\
  \vdots   \\
\sum_{i=1}^{\ell-1} \sum_{j=1}^{J_i}  I_{ij}(\theta,\lambda_{\ell-1}) z_{ij}
(\theta) \ge 1, \\
\sum_{i<\ell} \sum_{j=1}^{J_i} I_{ij}(\theta,\lambda_\ell) z_{ij}(\theta)
+ \sum_{j \notin J(\theta)} I_{\ell j}(\theta,\lambda_\ell)
z_{\ell j}(\theta) \ge 1.
\end{array} \right. \end{equation}
By Lemma 1, (\ref{mplemma1}) is true for all $\lambda\in\Lambda_{\ell}$. Therfore, 
$\liminf_{N\rightarrow\infty} R_{N}(\theta)/\log N
\ge \sup_{\lambda \in \Lambda_\ell}z(\theta,\ell,\lambda)$,
for all $\theta \in \Theta_\ell$. The proof is completed, if we can show that 
\begin{equation}\label{twompsame}
z(\theta,\ell) = \sup_{\lambda \in \Lambda_\ell}z(\theta,\ell,\lambda).
\end{equation}
If $Z=\{z_{ij}(\theta): j=1,\cdots,J_i~ {\rm for}~ i<\ell, {\rm and}~ 
j\not\in J(\theta),~ i=\ell\}$ satisfy (\ref{mpconstraint}), then
$Z$ also satisfy (\ref{mplemma1}). 
%\begin{equation}\label{compareconstraint}
%\inf_{\theta' \in \Theta_k} I_{ij}(\theta,\theta')
%\leq I_{ij} (\theta,\lambda_k),~ {\rm for}~ i\leq k;~~ 
%\inf_{\theta'\in B_{\ell}(\theta)} I_{ij}(\theta,\theta')
%\leq I_{ij}(\theta,\lambda_\ell),~ {\rm for}~ i\leq \ell.
%\end{equation}
Thus
\begin{equation}\label{glb}
 z(\theta,\ell) \geq \sup_{\lambda \in \Lambda_\ell}z(\theta,\ell,\lambda).
\end{equation}

Because $I_{ij}(\theta,\theta')$ are continuous with respect to $\theta'$,
the infimums in (\ref{mpconstraint}) are attained for some $\bar{\lambda}\in\bar{\Lambda}_{\ell}$, 
the closure of $\Lambda_{\ell}$.  Choose a sequence of 
$\lambda(n) =(\lambda_{1}(n),\cdots,$ $\lambda_\ell(n)) \in \Lambda_\ell$ such that it converges to the 
$\bar{\lambda}=(\bar{\lambda}_1,\cdots,\bar{\lambda_{\ell}})$. Note that $\bar{\lambda}$ depends on
some feasible $z$ satisfying (\ref{mpconstraint}).

Let $z_n =(z_{11}(n),\cdots,z_{\ell J_\ell}(n))$ be the solution of 
(\ref{mp}) satisfying (\ref{mplemma1}) with $\lambda=\lambda(n)$. Set
\[
c_{ij}(n)=\max\{I_{ij}(\theta,\lambda_1(n))/ I_{ij}(\theta,
\bar{\lambda}_1),\ldots,I_{ij}(\theta,\lambda_\ell(n))/ I_{ij}(\theta,
\bar{\lambda}_\ell) \}.
\]
By the continuity of $I_{ij}$, we have
\begin{equation}\label{limitcij}
\lim_{n\rightarrow\infty}c_{ij}(n)=1,~~~{\rm for}~1\le i \le \ell.
\end{equation}

In view of $\sum_{ij}c_{ij}(n)z_{ij}(n)I_{ij}(\theta,\bar{\lambda}_i)=\sum_{ij}z_{ij}(n)I_{ij}
(\theta,\lambda_i(n))$ for $i,j$ in an appropriate index set, we see that $\{c_{ij}(n)z_{ij}(n)\}$
satisfy (\ref{mpconstraint}). Hence,
\begin{eqnarray*}
&~& \Big[\max_{1\le i\le \ell,~1\le j\le J_i} c_{ij}(n) \Big] 
z(\theta,\ell,\lambda_n) \\
&\ge& \sum_{i<\ell} \sum_{j=1}^{J_i} [\mu^*(\theta)-\mu_{ij}(\theta)]
 c_{ij}(n)z_{ij}(n) + \sum_{j \notin J(\theta)}[\mu^*(\theta)-\mu_{\ell j}(\theta)]
 c_{\ell j}(n)z_{\ell j}(n) \ge z(\theta,\ell). 
\end{eqnarray*}
By (\ref{limitcij}), we have 
$\sup_{\lambda \in \Lambda_\ell} z(\theta,\ell,\lambda) \geq
z(\theta,\ell),$ which combined with (\ref{glb}) implies (\ref{twompsame}).

\subsection{Proof of Lemma 2}
By (\ref{localsup}), there exists
$\delta' > 0$ such that
\begin{equation}\label{localsup'}
E_{\pi_{ij}(\theta_0)} \Big[ \sup_{\tilde \theta \in N_{\delta'}(\theta')}
| \ell_{ij}(X_{ij0},X_{ij1};\theta',\wtd \theta) | \Big] < \varepsilon 
\end{equation}
for all $i,j$ and $\theta' \in \Theta$, $\varepsilon > 0$ to be specified later. Let
\begin{eqnarray}\label{xitilde}
\wtd \xi_{1jt} & = & \inf_{\lambda \in N_{\delta'}(\theta')} 
\ell_{1j}(X_{1j(t-1)},X_{1jt};\theta_0,\lambda) \cr
& = & \ell_{1j}(X_{1j(t-1)},X_{1jt};\theta_0,\theta') - 
\sup_{\lambda \in N_{\delta'}(\theta')} \ell_{1j}(X_{1j(t-1)},X_{1jt};\theta',\lambda).
\end{eqnarray}
Since $\eta:= \sum_{j=1}^{J_1} I_{1j}(\theta_0,\theta') > 0$, we 
can select $\delta' > 0$ to satisfy (\ref{localsup'}) with $\varepsilon < \eta/J_1$. 
 Then by (\ref{localsup'})-(\ref{xitilde}), it follows that 
\[
n_0^{-1} E_{\pi(\theta_0)} \Big[ \sum_{j=1}^{J_1} \sum_{t=1}^{n_0} 
\wtd \xi_{1jt} \Big]
\geq \sum_{j=1}^{J_1} I_{1j}(\theta_0,\theta') - J_1 \varepsilon \geq \eta - 
J_1 \varepsilon
> 0.
\]
By the Harris recurrence condition A3 and the law of large numbers, it follows that
\[
P(A) \rightarrow 1 \quad {\rm as} \ N \rightarrow \infty, \quad {\rm where} \
A =  \Big\{ \sum_{j=1}^{J_1} \sum_{t=1}^{n_0} \wtd 
\xi_{1jt} > 0 \Big\}.
\]
In the event $A$, the likelihood at $\theta_0$ is larger than all 
$\lambda \in N_{\delta'}(\theta')$ and hence (\ref{epmle}) holds. 

To prove (\ref{testingregretinf}), we extend (\ref{xiq}) and define
\begin{eqnarray}\label{xibreve}
\breve{\xi}_{ijt} & = & \inf_{\theta \in N_{\delta'}(\theta_0),\lambda \in
N_{\delta'}(\theta')} \log \Big[ {p_{ij}(X_{ij(t-1)},X_{ijt};\theta)
\over p_{ij}(X_{ij(t-1)},X_{ijt};\lambda)} \Big] \cr
& \geq & \ell_{ij}(X_{ij(t-1)}, X_{ijt};\theta_0,\theta')-
\sup_{\theta \in N_{\delta'}(\theta_0)}
| \ell_{ij}(X_{ij(t-1)},X_{ijt};\theta_0,\theta)| \cr
& & \qquad -
\sup_{\lambda \in N_{\delta'}(\theta')}
| \ell_{ij}(X_{ij(t-1)},X_{ijt};\theta',\lambda)|. 
\end{eqnarray}
Let $\theta' \in \bar \Theta_{kj_0}$ for some $k < \ell$. By A2, we can select
$0 < \varepsilon < I_{kj_0}(\theta_0,\theta')/2J_k$ and hence by (\ref{localsup'}) and 
(\ref{xibreve}), we have 
$E_{\pi(\theta_0)} \Big( \sum_{j=1}^{J_k} \breve{\xi}_{kjt} \Big) \geq I_{kj_0}
(\theta_0,\theta')-2 J_k \varepsilon > 0.$

By (\ref{teststat}) and (\ref{xibreve}), it follows that
\begin{equation}\label{teststatlemma2}
\inf_{\lambda \in N_{\delta'}(\theta')} \log U_k(\bn;\lambda) \geq 
\log F_k(N_{\delta'}(\theta_0)) +\sum_{i=1}^k 
\sum_{j=1}^{J_i} \log v_{ij} + 
\sum_{i=1}^k \sum_{j=1}^{J_i} \sum_{t=1}^{n_{ij}}\breve{\xi}_{ijt}, 
\end{equation}
where $v_{ij} = \inf_{x,\theta,\lambda} [\nu_{ij}(x;\theta)/\nu_{ij}
(x;\lambda)]$.  By (\ref{teststatlemma2}), $\tau_{kj}
\leq n_{kj}-m_{kj}$, where $\bn = (n_{ij})$ is the sample size needed for
$\sum_{i=1}^k \sum_{j=1}^{J_i} \sum_{t=1}^{n_{ij}}
\xi_{ijt}$ to cross the threshold $c := \log N - \sum_{i=1}^k \sum_{j=1}^{J_i} \log v_{ij}
-  \log F_k(N_{\delta'}(\theta_0))$ and $\bm = (m_{ij})$ is the sample size at the start of the
testing phase. Now follow arguments analogous to (\ref{logteststat}) - (\ref{tsincrement'}), 
we can prove the first half of (\ref{testingregretinf}).

Next, let us consider $k=\ell$. Let $f(\theta) = \mu_{\ell j_0}(\theta)-
\sup_{j \in J(\theta_0)} \mu_{\ell j}(\theta)$ for some $j_0 \not\in 
J(\theta_0)$. Then $f(\theta_0) < 0$. Conversely, $f(\theta') \geq 0$
for any $\theta' \in \Theta_{\ell j_0}$. By A1, $f$ is continuous with respect 
to $\theta$ and hence $\inf_{\theta' \in \Theta_{\ell j_0}} \| \theta_0 - \theta' \| > 0$.
The proof for second half of (\ref{testingregretinf}) then follows from the arguments similar
to those in the last two paragraphs of Section 6.1. 

\subsection{Extension of Wald's equation to Markovian rewards}

As we will be focusing on a single job $ij$ and
fixed parameters $\theta_0$, $\theta_q$ such that $\mu: = 
I_{ij}(\theta_0,\theta_q) > 0$
we will drop some of the references to $i$, $j$, $\theta_0$, $\theta_q$ and $q$ in
this subsection. This applies also to the notations in assumptions A3-A5.
Moreover, we shall use the notation $E(\cdot)$ as a short form of 
$E_{\theta_0}(\cdot)$ and $E_x(\cdot)$ as a short form of $E_{\theta_0}
(\cdot|X_0=x)$. 
Let $S_n = \xi_1 + \cdots + \xi_n$, where $\xi_k  =\log [p_{ij}(X_{k-1},
X_k;\theta_0)/p_{ij}(X_{k-1},$ $X_k;\theta_q)]$ has stationary mean under
$P_{\theta_0}$ and let $\tau$ be a stopping-time. We shall establish
Wald's equation
\begin{equation} \label{waldmarkovian}
E S_\tau = [\mu+o(1)] E \tau
\end{equation}
for Markovian rewards. 

By (\ref{hrecurrent}), we can augment the
Markov additive process and create a split chain containing an atom, 
so that increments in $S_n$ between visits to the atom are independent.
More specifically, we construct stopping-times $0 < \kappa(1) < \kappa(2)
< \cdots$ using an auxiliary randomization procedure such that
\begin{equation} \label{regeneration} \qquad
P \{ X_{n+1} \in A, \kappa(i)=n+1 | X_n=x, \kappa(i) > n \geq \kappa(i-1) \} = \cases{
\alpha \varphi(A) & if $x \in G$, \cr
0 & otherwise. }
\end{equation}
Then by Lemma 3.1 of Ney and Nummelin (1987),

\smallskip
(i) $\{ \kappa(i+1)-\kappa(i): i=1,2,\ldots \}$ are i.i.d. random variables.

\smallskip
(ii) the random blocks $\{ X_{\kappa(i)},\ldots,X_{\kappa(i+1)-1} \}$,
$i=1,2,\ldots,$ are independent and

\smallskip
(iii) $P \{ X_{\kappa(i)} \in A | {\cal F}_{\kappa(i)-1} \} = \varphi(A)$, where
${\cal F}_n$=$\sigma$-field generated by $\{ X_0,\ldots,X_n \}$.

\smallskip \noi Define
$\kappa=\kappa(1)$. By (ii)-(iii), $E_\varphi(S_\kappa-\kappa \mu)=0$. We preface the proof of (\ref{waldmarkovian}) with
the following preliminary lemmas, whose proofs are given in Chan, Fuh and
Hu (2005). 
  
\medskip
\begin{lemma}
Let $\gamma(x) = E_x (S_\kappa-\kappa \mu)$. Then
$Z_n = (S_n-n \mu)+\gamma(X_n)$ is a martingale with respect to ${\cal F}_n$. 
Hence 
\begin{equation} \label{markovidentity}
E S_{\tau} = \mu (E \tau) - E[\gamma(X_\tau)] + E[\gamma(X_0)].
\end{equation}
\end{lemma}
\begin{lemma}
Under {\rm A3-A5},
$$
|\gamma(x)| \leq \bar{b}^{-1}[V(x)+b+(V^*+b)V^*(\alpha^{-1}+1)](K+1+|\mu|),
$$
where $\alpha$ satisfies {\rm (\ref{hrecurrent})}, $V^*$ is defined in {\rm (\ref{vij})}
and $K$ is defined in {\rm (\ref{a51})}. 
\end{lemma}

\medskip
Let $W_i = |\gamma(X_{\kappa(i)})| + \cdots + |\gamma(X_{\kappa(i+1)-1})|$, for
$i \geq 1$. Then by A3-A5, Lemma 4 and its proof, and (i)-(iii), $W_1,W_2,\ldots$
are i.i.d. with finite mean while by (\ref{vij}),
$W_0 := | \gamma(X_0)|+\cdots+|\gamma(X_{\kappa(1)-1})|$ also has finite mean.

\medskip
\begin{lemma}
Let $M_n = \max_{1 \leq k \leq n} W_k$. Then for any stopping-time $\tau$,
$E(M_\tau) = o(E \tau)$. 
\end{lemma}

\medskip \noi {\bf Proof of (\ref{waldmarkovian})}. By Lemma 5,
$E|\gamma(X_\tau)| + E|\gamma(X_0)| = o(E \tau)$, and  (\ref{waldmarkovian})
follows from (\ref{markovidentity}). $\Box$

%----------  END DOCUMENT BODY -------- BEGIN BIBLIOGRAPHY ----------

\end{document}